\input amstex
\magnification=\magstep1
\input epsf
\baselineskip=13pt
\documentstyle{amsppt}
\vsize=8.7truein
\CenteredTagsOnSplits
\NoRunningHeads
\def\today{\ifcase\month\or
  January\or February\or March\or April\or May\or June\or
  July\or August\or September\or October\or November\or December\fi
  \space\number\day, \number\year}
\def\dist{\operatorname{dist}}
\def\conv{\operatorname{conv}}
\def\spa{\operatorname{span}}
\def\PP{{\bold P}}

\def\Mat{\operatorname{Mat}}
\topmatter
\title The Distribution of Values in the Quadratic Assignment Problem \endtitle
\author Alexander Barvinok and Tamon Stephen\endauthor
\address Department of Mathematics, University of Michigan, Ann Arbor,
MI 48109-1109 \endaddress
\email barvinok$\@$umich.edu, tamon$\@$umich.edu  \endemail
\date  \today \enddate
\thanks This research was partially supported by NSF Grant DMS 9734138.
\endthanks
\abstract We obtain a number of results regarding the distribution 
of values of a quadratic function $f$ on the set of $n \times n$ 
permutation matrices (identified with the symmetric group $S_n$) around its 
optimum (minimum or maximum). In particular, we estimate the 
fraction of permutations $\sigma$ such that $f(\sigma)$ lies within a 
given neighborhood of the optimal value of $f$.
We identify some ``extreme'' functions $f$ (there are 4
of those for $n$ even and 5 for $n$ odd) such that the distribution of
every quadratic function around its optimum is a certain ``mixture'' 
of the distributions of the extremes and describe 
a natural class of functions (which includes, 
for example, the objective function in the Traveling Salesman Problem)
 with a relative abundance of 
near-optimal permutations. In particular, we identify a 
large class of functions $f$ with the property that permutations in the 
vicinity of the optimal permutation 
(in the Hamming metric of $S_n$) tend to produce 
near optimal values of $f$ (such is, for example,
the objective function in the symmetric Traveling 
Salesman Problem) and show that for general $f$, just the opposite behavior 
may take place: an average permutation in the vicinity of the 
optimal permutation may be much worse than an average permutation in 
the whole group $S_n$.  
\endabstract
\endtopmatter
\document

\head 1. Introduction \endhead

The Quadratic Assignment Problem (QAP for short) is an optimization
problem on the symmetric group $S_n$ of $n!$ permutations of an 
$n$-element set.
The QAP is one of the hardest problems of combinatorial optimization,
whose special cases include the Traveling Salesman Problem (TSP) among 
other interesting problems. 

Recently the QAP has been of interest to many people.
An excellent survey of recent results is found in [5].
Despite this work, it is still extremely difficult to solve
QAP's of size $n=20$ to optimality, and the solution to a
QAP of size $n=30$ is considered noteworthy, see, for example, 
[1] and [4]. Moreover, it appears that essentially no positive 
approximability results for the general QAP are known, although 
some ``bad news'' (non-approximability) and approximability for 
special classes have been established, see [3] and [2]. 

The goal of this paper is to study the distribution of values of the 
objective function of the QAP. 
We hope that our results would allow one on one hand to understand
the behavior of the local search heuristic, and, on the other hand, to get 
guaranteed approximations to the optimum using some simple algorithms 
based on random or partial enumeration with guaranteed complexity
bounds. In particular, we estimate how well the sample optimum 
from a random sample of a given size approximates the global 
optimum.
\subhead (1.1) The Quadratic Assignment Problem \endsubhead
Let $\Mat_n$ be the vector space of all real $n \times n$ matrices
$A=(a_{ij})$, $1 \leq i, j \leq n$ and let $S_n$ be 
the set of all permutations $\sigma$ of the set $\{1, \ldots, n\}$.
There is an action of $S_n$ on the space $\Mat_n$ by simultaneous 
permutations of rows and columns: we let $\sigma(A)=B$, where 
$A=(a_{ij})$ and $B=(b_{ij})$, provided $b_{\sigma(i)\sigma(j)}=a_{ij}$
for all $i,j=1,\ldots, n$. One can check that $(\sigma \tau)A=\sigma (\tau A)$
for any two permutations $\sigma$ and $\tau$.
There is a standard scalar product on $\Mat_n$:
$$\langle A, B \rangle=\sum_{i,j=1}^n a_{ij} b_{ij} \quad 
\text{where} \qquad A=(a_{ij}) \quad \text{and} \quad
B=(b_{ij}).$$  
Let us fix two matrices $A=(a_{ij})$ and $B=(b_{ij})$
and let 
us consider a real-valued function $f: S_n \longrightarrow {\Bbb R}$ 
defined by 
$$f(\sigma)=\langle B, \sigma(A) \rangle =\sum_{i,j=1}^n 
b_{\sigma(i) \sigma(j)} a_{ij}=\sum_{i,j=1}^n b_{ij} 
a_{\sigma^{-1}(i) \sigma^{-1}(j)}  \tag1.1.1$$
The problem of finding a permutation $\sigma$ where the maximum or 
minimum value of $f$ is attained is known as the {\it Quadratic Assignment
Problem}. It is one of the hardest problems of Combinatorial Optimization. 
From now on we assume that $n \geq 4$.
\bigskip
 In this paper, we study the distribution 
of values of $f$ from the optimization perspective: 
\bigskip
$\bullet$ How ``steep'' or 
how ``flat'' can the optimum of $f$ be?
\medskip
$\bullet$ How many values of $f$ lie within a given distance to the optimum?  
\medskip
$\bullet$ When can we hope to improve the value of $f(\sigma)$ by modifying 
$\sigma$ slightly?
\bigskip
To formulate the questions rigorously (and to answer them), we introduce 
the standard Hamming metric on the symmetric group $S_n$.
\definition{(1.2) Definitions} For two permutations $\tau, \sigma \in S_n$,
let the distance $\dist(\sigma, \tau)$ be the number of indices 
$1\leq i \leq n$ where $\sigma$ and $\tau$ disagree:
$$\dist(\tau, \sigma)=|i: \sigma(i) \ne \tau(i)|.$$
One can observe that the distance is invariant under the left and 
right actions of $S_n$:
$$\dist(\sigma \sigma_1, \sigma \sigma_2)=\dist(\sigma_1, \sigma_2)=
\dist(\sigma_1 \sigma, \sigma_2 \sigma)$$
for all $\sigma_1, \sigma_2, \sigma \in S_n$.

For a permutation $\tau$ and an integer $k>1$, we consider the 
``$k$-th ring'' around $\tau$:
$$U(\tau, k)=\bigl\{ \sigma \in S_n: \dist(\sigma, \tau)=k \bigr\}.$$  
\enddefinition
In particular, we are interested in the distribution of values of $f$ 
in the set $U(\tau, k)$, where $\tau$ is an optimal permutation.
\subhead (1.3) The generalized problem \endsubhead
Our approach produces essentially identical results for a more 
general problem, where we are given a 4-dimensional array 
$C=\Bigl\{c^{ij}_{kl}: 1 \leq i,j, k,l \leq n\Bigr\}$ 
of $n^4$ real numbers and the 
function $f$ is defined by
$$f(\sigma)=\sum_{i,j=1}^n c^{ij}_{\sigma(i) \sigma(j)}. 
\tag1.3.1$$ 
If $c^{ij}_{kl}=a_{ij} b_{kl}$ for some matrices $A=(a_{ij})$ and
$B=(b_{kl})$, in which case we write $C=A \otimes B$, we get the 
special case (1.1.1) we started with.
\bigskip
The main idea of our approach is as follows. Let 
$$\overline{f}={1 \over n!} \sum_{\sigma \in S_n} f(\sigma)$$
be the average value of $f$ on the symmetric group and let 
$f_0=f-\overline{f}$. Hence the average value of $f_0$ is 0 and we study 
the distribution of values of $f_0$ around its maximum (the problem 
with minimum instead of maximum is completely similar).
Now, as long as the 
distribution of values of $f_0$ is concerned, without loss of 
generality we may assume 
that $f_0$ attains its maximum on the identity permutation $e$, so 
that $f_0(e) \geq f_0(\sigma)$ for all $\sigma \in S_n$. 
Let us define a function $g: S_n \longrightarrow {\Bbb R}$, 
which we call the {\it central projection} (with the term coming 
from the representation theory) of $f$ by
$$g(\sigma)={1 \over n!} \sum_{\omega \in S_n} f_0(\omega^{-1} \sigma \omega). 
\tag1.4$$ 
It turns out that $g$ attains its maximum on the identity permutation, 
that the average value of $g$ on $S_n$ is 0 and, moreover, the 
average values of $f_0$ and $g$ on the $k$-th ring $U(e, k)$ coincide for 
all $k$. In short, $g$ captures some important information about the 
distribution of $f$. 
The set of all functions $g$ obtained by central projection (1.4)
from all functions $f_0$ having maximum at the 
identity forms a 3-dimensional convex polyhedral cone.
We describe this cone, identifying its extreme rays 
(there are $4$ for even $n$ and $5$ for odd $n$), which provide 
us with some ``extreme'' types of distribution.
Hence we study the distribution of values of $g$, which is a much easier 
problem. Once 
the distribution of values of $g$ is understood,
using (1.4), we infer various facts about the distribution of values of $f$.  

We remark that it is easy to compute 
the average value $\overline{f}$ of $f$ given by (1.1.1) or by (1.3.1).
\proclaim{(1.5) Lemma} Let $f: S_n \longrightarrow {\Bbb R}$ be 
a function defined by 
$$f(\sigma)=\langle B, \sigma(A) \rangle$$ 
for some matrices $A=(a_{ij})$ and $B=(b_{ij})$. Let 
$$\overline{f}={1 \over n!} \sum_{\sigma \in S_n} f(\sigma)$$
be the average value of $f$ on the symmetric group $S_n$.
Let us define 
$$\split &\alpha_1=\sum_{1 \leq i \ne j \leq n} a_{ij}, \quad
 \alpha_2=\sum_{i=1}^n a_{ii}\quad \text{and} \\ 
&\beta_1=\sum_{1\leq i \ne j \leq n} b_{ij}, \quad
 \beta_2=\sum_{i=1}^n b_{ii}. \endsplit$$
Then 
$$\overline{f} = {\alpha_1 \beta_1 \over n(n-1)} + 
{\alpha_2 \beta_2 \over n}.$$
\endproclaim
Similarly, if $f$ is a function (1.3.1) of the generalized problem,
then 
$$\overline{f}={1 \over n(n-1)}\sum_{1 \leq i \ne j \leq n} 
\sum_{1 \leq k \ne l \leq n} c^{ij}_{kl} +{1 \over n} \sum_{1 \leq i, l \leq n}
c^{ii}_{ll}.$$

We prove Lemma 1.5 in Section 6.
\subhead (1.6) Notation \endsubhead We often denote by $c$ some positive 
constant whose precise value is not of particular importance to us.
If $F$ and $G$ are non-negative functions of a positive integer $n$, we 
write $F=O(G)$ if $F(n) \leq cG(n)$ for some $c>0$ and all sufficiently 
large $n$. Similarly, we write $F(n)=\Omega(G)$ if $F(n) \geq cG(n)$ 
for some constant $c>0$ and all sufficiently large $n$. We denote by 
$e$ the identity permutation in $S_n$. 
We denote by $|X|$ the cardinality of a finite set $X$ and by $\conv A$ the 
convex hull of the set $A$ in Euclidean space.
Given a function 
$f: S_n \longrightarrow {\Bbb R}$, we denote by $\overline{f}$ its 
average value on $S_n$:
$$\overline{f}={1 \over n!} \sum_{\sigma \in S_n} f(\sigma)$$
and by 
$$f_0=f-\overline{f}$$ 
the ``shifted'' function with 0 average.
Our results concern the function $f_0$.
\bigskip
The paper is organized as follows.
In Sections 2-5, we state our results about the number of near-optimal 
permutations. In Sections 6-11, we prove those results and describe 
certain ``extreme'' distributions.
 We give an informal preview 
of our results below. In what follows, $\tau$ is an optimal permutation such 
that $f_0(\tau) \geq f_0(\sigma)$ for all $\sigma \in S_n$. Since 
we consider the shifted function, the minimization and maximization 
problems are completely similar.  

In Section 2, we consider a special case of the problem where matrix 
$A$ is symmetric, has constant row and column sums and a constant diagonal
(of course, $A$ and $B$ are interchangeable).
For example, the symmetric TSP belongs to this class. 
The interesting feature of 
this special case is what we call the ``bullseye'' distribution of values 
of $f_0$ around its maximum. It turns out that the average 
value of $f_0$ over the $k$-th ring $U(\tau, k)$ (see Definitions 1.2)
around an optimal permutation $\tau$ steadily improves as the ring 
contracts to $\tau$. The proof is given in Section 8.
This is also the simplest case to analyze. It turns out that the set 
of all possible central projections $g$ (see (1.4)) is one-dimensional.

In Section 3, we consider a more general case of a not necessarily 
symmetric matrix $A$ with constant row and column sums and 
a constant diagonal. For example, 
the asymmetric TSP belongs to this class.
We call this case ``pure'' since the objective 
function $f$ lacks the component that can be attributed to the 
Linear Assignment Problem. Although we don't have the bullseye 
distribution of Section 2, we can
provide some guarantees for the number of reasonably
good permutations $\sigma$. Thus, for any $\alpha>1$ the probability that 
a random permutation $\sigma \in S_n$ satisfies 
$\displaystyle f_0(\sigma) \geq {\alpha \over n^2} f_0(\tau)$ is at least 
$\Omega(n^{-2})$. Furthermore,
for any $\epsilon>0$ the probability that a random permutation $\sigma$ 
satisfies $\displaystyle f_0(\sigma) \geq n^{-\epsilon} f_0(\tau)$ is 
``mildly exponential'', that is at least of the order of $\exp\{-n^{c}\}$
for some constant $c=c(\epsilon)<1$.
The proof is given in Section 9. It turns out that the set 
of all central projections $g$, maximized at the identity,
forms a 2-dimensional cone. The  
extreme rays provide us with the extreme types of distributions, which,
although not as good as the ``bullseye'' distribution of Section 2,
still quite reasonable, especially 
compared with types of distributions we encounter in general symmetric 
QAP.  

In Section 4, we consider the symmetric Quadratic Assignment Problem, where 
matrix $A$ (or, equivalently $B$) is symmetric. This case turns out to be 
very different in many respects from the special cases of Sections 2 and 3. 
It turns out that the ``bullseye'' distribution is no longer the law. 
We present a simple example of function $f_0$ where the average value 
of $f_0$ over the $k$-th ring $U(\tau, k)$ of an optimal permutation $\tau$ 
is much worse than the average over the whole group $S_n$ even for
small $k$. We call such a distribution a ``spike''. 
We argue that at least for the generalized 
problem (1.3), the number of near-optimal permutations is much smaller than
in the pure case of Section 3.
The proofs are given in Section 10. It turns out that the set of 
all central projections (1.4) forms a 2-dimensional cone whose 
extreme rays provide us with the extreme types of distributions. 
One of those rays turns out to have an extreme ``spike'' distribution.

In Section 5, we consider the general Quadratic Assignment Problem.
As in Section 3, we
prove that for any $\alpha>1$ the probability that 
a random permutation $\sigma \in S_n$ satisfies 
$\displaystyle f_0(\sigma) \geq {\alpha \over n^2} f_0(\tau)$ is at least 
$\Omega(n^{-2})$, although with a worse constant than in Section 
3. 
We prove that for any $\epsilon>0$ there is a constant $c(\epsilon)<1$ 
such that 
the probability that a random permutation $\sigma$ 
satisfies $\displaystyle f_0(\sigma) \geq n^{-1-\epsilon} f_0(\tau)$ 
is at least of the order $\exp\{-n^{c}\}$ (mildly exponential).
The proofs are given in Section 11. It turns out that the set 
of central projections (1.4) forms a 3-dimensional polyhedral cone 
with 4 extreme rays when $n$ is even and 5 extreme rays when $n$ is odd.
In a sense, those extreme rays describe all ``extreme'' distributions 
that one may encounter in the general Quadratic Assignment Problem.

In Section 6, we prove some preliminary technical results. In Section 7,
we review the necessary facts from the representation theory of 
the symmetric group, which we use essentially in our approach.

\head 2. The Bullseye Case \endhead

Our analysis
of the Quadratic Assignment Problem is the simplest 
in the following special case (it also exhibits some features 
absent in the general case).
Suppose that the matrix $A=(a_{ij})$ is symmetric and has constant row 
and column
sums and a constant diagonal:
$$\split a_{ij}=a_{ji} \quad &\text{for all} \quad 1\leq i,j \leq n;  \\
\\ &\text{ } \\
 &\text{for some} \quad a  \\ 
\sum_{i=1}^n a_{ij}=a
\quad &\text{for all} \quad j=1, \ldots, n \quad \text{and} \\
 \sum_{j=1}^n a_{ij}=a 
\quad &\text{for all} \quad i=1, \ldots, n; \\
\\ &\text{ } \\
a_{ii}=b  \quad &\text{for some} \quad b \quad \text{and all} \quad 
i=1, \ldots, n. \endsplit$$
For example, 
$$A=\left( \matrix 0 & 1 & 0 & \ldots & 0 & 1 \\
                   1 & 0 & 1 & \ldots & 0 & 0 \\
                   0 & 1 & 0 & \ldots & 0 & 0 \\
                   \ldots & \ldots & \ldots & \ldots & \ldots \\
                   0 & 0 & 0 & \ldots & 0 & 1 \\
                   1 & 0 & 0 & \ldots & 1 & 0 \endmatrix \right), \quad
a_{ij}=\cases 1 &\text{if \ }|i-j|=1 \mod n \\ 0 &\text{otherwise}
\endcases$$
satisfies these properties and the corresponding optimization problem is the 
{\it Symmetric Traveling Salesman Problem}.
It turns out that the optimum has a characteristic ``bullseye'' feature
in the Hamming metric on $S_n$ (see Definition 1.2).
\proclaim{(2.1) Theorem} Suppose that the matrix $A$ is symmetric
and has constant row and column sums and a constant diagonal.
Let $f: S_n \longrightarrow {\Bbb R}$ be the function defined 
by (1.1.1) for $A$ and some matrix $B$. Let $\overline{f}$ be the average 
value of $f$ on $S_n$, let $f_0=f-\overline{f}$ and
let $\tau \in S_n$ be an optimal permutation: $f_0(\tau)=\max_{\sigma \in S_n}
f_0(\sigma)$. 
For $k \geq 0$ let 
$$U(\tau; k)=\bigl\{\sigma: \dist(\sigma, \tau)=k \bigr\}$$ be the 
$k$-th ``ring'' around $\tau$ and let 
$$\alpha(n,k)={(n-k)^2-3(n-k) \over n^2-3n}.$$
Then
$${1 \over |U(\tau; k)|} \sum_{\sigma \in U(\tau; k)} 
f_0(\sigma) 
 \geq \alpha(n,k) f_0(\tau).$$  
\endproclaim
We prove Theorem 2.1 in Section 8. 
\subhead (2.2) The ``bullseye'' distribution. 
Connections with the local search \endsubhead
It follows from our proof that 
we have almost equality in the formula of Theorem 2.1. We observe that as 
the ring $U(\tau; k)$ contracts to the optimal permutation $\tau$, 
the {\it average} value of $f$ on the ring steadily improves.
$$\epsffile{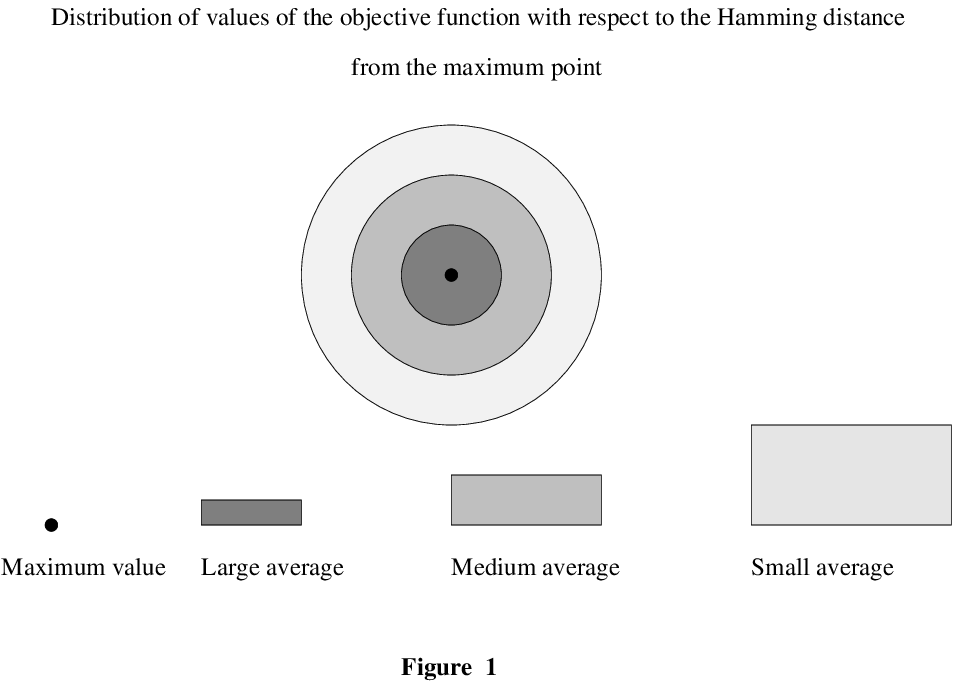}$$
It is easy to construct examples where {\it some} values
of $f$ in a very small neighborhood of the optimum are particularly bad,
but as follows from Theorem 2.1, such values are
relatively rare. In our opinion,
this provides some justification for the local search heuristic, 
where one starts from a permutation and tries to 
improve the value of the objective function by searching a small
neighborhood of the current solution. Indeed, if we had the value 
of $f_0(\sigma)$ for {\it each} $\sigma  \in U(\tau, k)$ equal to 
$\alpha(n,k)f_0(\tau)$, then the local search 
would have converged to the optimum in $O(n)$ steps, since each step 
would have brought us to a smaller neighborhood of the optimal solution.
Instead, we have that the average value over $U(\tau, k)$ is (almost) 
equal to $\alpha(n,k)f_0(\tau)$. We can no 
longer guarantee that the local search converges fast (or even 
converges) to the optimal 
solution (after all, our problem includes the Traveling Salesman Problem 
as a special case and hence is NP-hard), but it plausible that the 
local search behaves reasonably well for an ``average'' optimization 
problem. This agrees with the empirical evidence that the local search works 
well for the Traveling Salesman Problem.

Incidentally, one can prove that 
the same type of the ``bullseye'' behavior is observed for the
 Linear Assignment Problem and some other polynomially solvable problems,
such as the weighted Matching Problem.
\bigskip
Estimating the size of the ring $U(\tau, k)$, we get the following
result.
\proclaim{(2.3) Theorem} 
Suppose that the matrix $A$ is symmetric and 
has constant row and column sums and a constant diagonal. Let 
$f: S_n \longrightarrow {\Bbb R}$ be the function defined by (1.1.1) for 
$A$ and some matrix $B$, let $\overline{f}$ be the 
average value of $f$ on $S_n$ and let $f_0=f-\overline{f}$.
 Let $\tau$ be an optimal permutation:
$f_0(\tau)=\max_{\sigma \in S_n} f_0(\sigma)$. 
Let us choose an integer $3 \leq k \leq n-5$ and a number $0<\gamma<1$ and let 
$$\beta(n,k)={k^2-3k \over n^2-3n}.$$
The probability that a random permutation $\sigma \in S_n$ 
satisfies the inequality 
$$f_0(\sigma) \geq \gamma \beta(n,k)f_0(\tau)$$
is at least
$${(1-\gamma)\beta(n,k) \over 3k!}.$$
\endproclaim
We prove Theorem 2.3 in Section 8.

Our results can be generalized in a quite 
straightforward way to functions $f$ defined by (1.3.1), if 
we assume that for any $k$ and $l$ the matrix $A=(a_{ij})$, where 
$a_{ij}=c^{ij}_{kl}$, is 
symmetric with constant row and column sums and has a constant 
diagonal.

\head 3. The Pure Case \endhead 

In this Section, we consider a more general case of 
a not necessarily symmetric matrix $A$ having constant row and column 
sums and a constant diagonal:
$$\split &\text{ } \\
 &\text{for some} \quad a  \\ 
\sum_{i=1}^n a_{ij}=a
\quad &\text{for all} \quad j=1, \ldots, n \quad \text{and} \\
 \sum_{j=1}^n a_{ij}=a 
\quad &\text{for all} \quad i=1, \ldots, n; \\
\\ &\text{ } \\
a_{ii}=b  \quad &\text{for some} \quad b \quad \text{and all} \quad 
i=1, \ldots, n. \endsplit$$
For example, matrix 
$$A=\left( \matrix 0 & 1 & 0 & \ldots & \ldots & 0 \\
           0 & 0 & 1 & 0 &\ldots & 0 \\
           \ldots & \ldots &\ldots &\ldots &\ldots &\ldots \\
           0 & 0 & \ldots & \ldots & 0 & 1 \\
           1 & 0 &\ldots &\ldots & \ldots &0 \endmatrix \right).
\qquad a_{ij}=\cases 1 &\text{if \ }j=i+1 \mod n \\
                     0 & \text{otherwise} \endcases$$  
satisfies these properties and the corresponding optimization problem 
is the {\it Asymmetric Traveling Salesman Problem}.

We call this case pure, because as we remark in Sections 7 and 9, 
the objective function $f$ lacks the component attributed to the 
Linear Assignment Problem. More generally, an arbitrary objective 
function $f$ in the Quadratic Assignment Problem can be represented 
as a sum $f=f_1+f_2$, where $f_1$ is the objective function in a 
Linear Assignment Problem and $f_2$ is the objective function in some 
pure case. 

In this case we can no longer claim the bullseye distribution of 
Section 2 (the reasons are explained in Section 9), the distribution in this
case is not as bad as, for example, in the general symmetric QAP (see 
Section 4) and the estimates of the number of relatively 
good values we are able to prove are almost as good as 
those of Section 2.
\proclaim{(3.1) Theorem} Suppose that the matrix $A$ has constant row 
and column sums and a constant diagonal. Let $f: S_n \longrightarrow {\Bbb R}$ 
be the function defined by (1.1.1) for $A$ and some matrix $B$, let 
$\overline{f}$ be the average value of $f$ on $S_n$ and let 
$f_0=f-\overline{f}$. Let $\tau$ be an optimal permutation, so
$f_0(\tau)=\max_{\sigma \in S_n} f_0(\sigma)$.
Let us choose an integer $3 \leq k \leq n-5$ and a number $0<\gamma<1$ and let
$$\beta(n,k)={k^2-3k+1 \over n^2-3n+1}.$$
The probability that a random permutation $\sigma \in S_n$ satisfies 
the inequality 
$$f_0(\sigma) \geq \gamma \beta(n,k) f_0(\tau)$$
is at least
$${(1-\gamma) \beta(n,k) \over 10k!}.$$ 
\endproclaim 
In particular, by choosing an appropriate $k$, we obtain the following 
corollary.
\proclaim{(3.2) Corollary} 

\roster
\item Let us fix any $\alpha>1$. Then there exists a $\delta=\delta(\alpha)>0$ 
such that for all sufficiently large $n \geq N(\alpha)$ the probability that a 
random permutation $\sigma$ in $S_n$ satisfies the inequality
$$f_0(\sigma)\geq {\alpha \over n^2} f_0(\tau)$$
is at least $\delta n^{-2}$. In particular, one can choose 
$\displaystyle\delta=\exp\bigl\{-c \sqrt{\alpha} \ln \alpha \bigr\}$ for 
some absolute constant $c>0$.
\item Let us fix any $\epsilon>0$. Then there exists a 
$\delta=\delta(\epsilon)<1$ 
such that for all sufficiently large $n \geq N(\alpha)$ the probability that 
a random permutation $\sigma$ in $S_n$ satisfies the inequality 
$$f_0(\sigma) \geq n^{-\epsilon} f_0(\tau)$$
is at least $\exp\{-n^{\delta}\}$. In particular, one can choose 
any $\delta>1-\epsilon/2$.
\endroster
\endproclaim
{\hfill \hfill \hfill} \qed

We prove Theorem 3.1 in Section 9.

From Corollary 3.2, it follows that to get a permutation $\sigma$
which satisfies (1) for any fixed $\alpha$, we can use the 
following straightforward randomized algorithm: 
sample $O(n^2)$ random permutations
$\sigma \in S_n$, compute the value of $f$ and choose the best permutation. 
With the probability which tends to 1 as $n \longrightarrow +\infty$, 
we will hit the right permutation. The complexity of the algorithm is 
quadratic in $n$ for any $\alpha$, but the coefficient of $n^2$ grows 
as $\alpha$ grows. If we are willing to settle for 
an algorithm of a mildly exponential complexity 
of the type $\exp\{n^{\beta}\}$ for some $\beta<1$
we can achieve a better approximation (2) by searching through the 
set of randomly selected $\exp\{n^{\beta}\}$ permutations. We remark 
that no algorithm solving the Quadratic Assignment Problem (even 
in the special case considered in this section) with an exponential 
in $n$ complexity $\exp\{O(n)\}$ is known, although there is a 
dynamic programming algorithm solving the Traveling Salesman Problem in 
$\exp\{O(n)\}$ time.

Again, our results can be generalized in a quite 
straightforward way to functions $f$ defined by (1.3.1), if 
we assume that for any $k$ and $l$ the matrix $A=(a_{ij})$, where 
$a_{ij}=c^{ij}_{kl}$ has constant row and column sums and has a constant 
diagonal.

\head 4. The Symmetric Case \endhead

In this section, we assume that the matrix $A=(a_{ij})$ is symmetric,
that is 
$$a_{ij}=a_{ji} \quad \text{for all} \quad 1\leq i,j \leq n.$$
Overall, the distribution of values of $f$ 
turns out to be much more complicated when in the special cases described 
in Sections 2 and 3. First, we observe that generally one can not hope 
for the ``bullseye'' feature described in Section 2.2.
\subhead (4.1) The ``spike'' distribution \endsubhead
Let us choose an $n \times n$ matrix $A=(a_{ij})$, where 
$$a_{ij}=\cases 1 &\text{if\ } (ij)=(12) \quad \text{or \ } (ij)=(21) \\ 
0 &\text{otherwise,} \endcases$$
so 
$$A=\left( \matrix 0 & 1 & 0 & \ldots & 0 \\ 1 & 0 & 0 &\ldots & 0 \\
0 & 0 & 0 & \ldots & 0 \\ \ldots & \ldots & \ldots & \ldots & \ldots \\ 
0 & 0 & 0 & \ldots & 0 \endmatrix \right).$$
Let 
$$\gamma={-n^2+5n-8 \over 8(n-2)}$$ and let $B=(b_{ij})$, where 
$$b_{ij}=\cases 0 &\text{if\ } i=j \\ \gamma &\text{if \ } i \leq 2 
\text{ \ and \ }j \geq 3 \quad \text{or} \quad \text{if \ } j \leq 2 
\text{ \  and \ }i \geq 3 \\ 1/2 &\text{otherwise,} \endcases$$
so 
$$B=\left(\matrix 0 & 1/2 & \gamma  &\gamma 
&\ldots & \gamma \\ 1/2 & 0 & \gamma & \gamma 
&\ldots & \gamma \\ \gamma & \gamma & 0 & 1/2 & \ldots & 1/2 \\
\ldots &\ldots &\ldots &\ldots &\ldots &\ldots \\
\gamma & \gamma & 1/2 & 0 & \ldots & 1/2 \\
\\ \gamma &\gamma &1/2 
&\ldots &1/2  & 0 \endmatrix \right).$$
Let $f: S_n \longrightarrow {\Bbb R}$ be the function defined by (1.1.1). 
In Section 10, we prove the following properties of $f$.
\bigskip
$\bullet$ We have $\overline{f}=0$ for the average value of $f$ on $S_n$;
\medskip
$\bullet$ 
The maximum value of $f$ on $S_n$ is 1 and is attained, in particular, on the 
identity permutation $e$;
\medskip
$\bullet$ For the $k$-th ring $U(e,k)$ 
centered at the identity permutation $e$, 
we have 
$${1 \over |U(e,k)|} \sum_{\sigma \in U(e,k)} f(\sigma) 
\leq {-nk+k^2+3n-k-4 \over 2n-4}.$$
\bigskip
We observe that already for $k=4$ (a more careful analysis 
yeilds $k=3$) the average 
value of $f$ over $U(e,k)$ is negative for all sufficiently 
large $n$. 
Thus an average permutation in $U(e, 4)$ presents us 
with a choice worse than an average permutation in $S_n$. 
The distribution of values of $f$ turns out to be of 
the opposite nature to the bullseye distribution of Figure 1.
We call it the ``spike'' distribution.
$$\epsffile{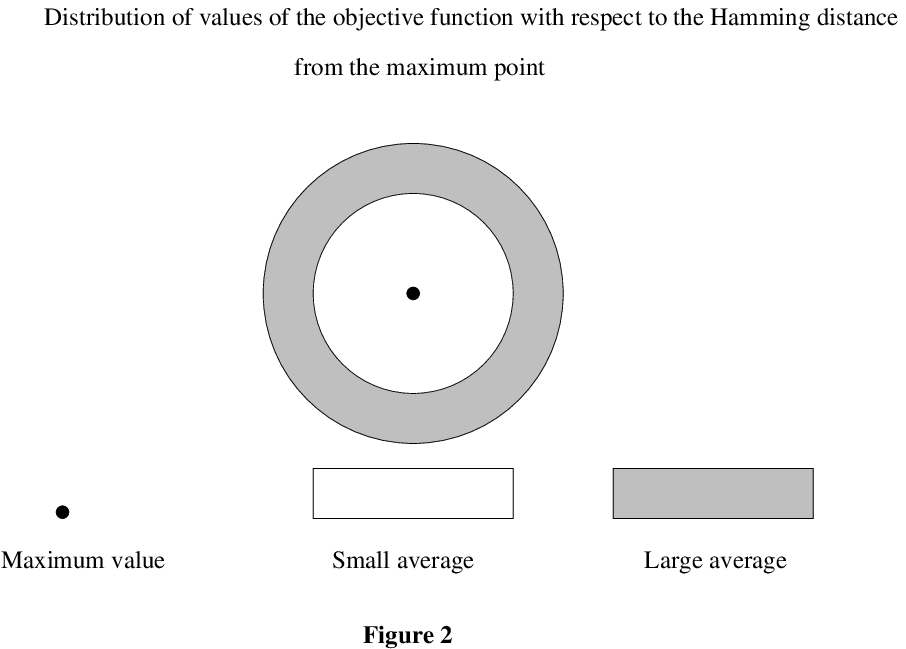}$$
Of course, in this particular case the optimization 
problem is very easy since the function $f$ attains only 
two different values. However, this may serve as an indication that 
complicated distributions are indeed possible and the local search 
may not work well for a general symmetric QAP.
Indeed, this is the case if we allow generalized
functions (1.3.1). 

In Section 10, we show that there exists 
a tensor $c^{ij}_{kl}$ with the property that $c^{ij}_{kl}=c^{ji}_{kl}$ 
for all $k$ and $l$ and all $i$ and $j$ such that for the corresponding 
function $f$ defined
by (1.3.1), 
we have 
$$f(\sigma)={-n\bigl(p(\sigma)-1)+p(\sigma)(p(\sigma)+1)+2t(\sigma)-4
\over 2n-4}, \tag4.1.1$$
where $p(\sigma)=\big|\{i: \sigma(i)=i \}\big|$ is the number of fixed 
points of the permutation and $t(\sigma)=\big|\bigl\{i<j: \sigma(i)=j$ 
and $\sigma(j)=i \bigr\}\big|$ is the number of 2-cycles in the permutation.
We show that $\overline{f}=0$ and that the maximum value 1 of $f$ is 
attained at the identity permutation $e$ (where $p(\sigma)=n$ and
$t(\sigma)=0$) and, for even $n$, 
on the permutations 
that consist of $n/2$ transpositions (where $p(\sigma)=0$ and 
$t(\sigma)=n/2$). On the other hand, for any fixed 
$n-3 \geq k \geq 3$ and all $n \geq 5$, 
the value of $f(\sigma)$ with $\sigma \in U(e, k)$ is negative.
\subhead (4.2) Scarcity of relatively good values \endsubhead
Unfortunately, we are unable to present an example of the symmetric QAP
which beats the bound of Theorem 3.1 but we can construct such
an example for the generalized problem (1.3). In Section 10, we prove
that for any $3 \leq m \leq n$, there exists a tensor $c^{ij}_{kl}$ such 
that $c^{ij}_{kl}=c^{ji}_{kl}$ for all $k$ and $l$ and such that 
for the corresponding function $f$ we have
$$f(\sigma)={p^2(\sigma)-mp(\sigma)+2t(\sigma)+m-3 \over n^2 -nm +m-3},
\tag4.2.1$$
where $p(\sigma)$ is the number of fixed points in $\sigma$ and 
$t(\sigma)$ is the number of 2-cycles in $\sigma$.
We show that $\overline{f}=0$ and that $f(e)=1$ is the maximum value of 
$f$. 

Let us fix any $0<\delta<1$ and let us choose some $m$ such that
$n^{1-\epsilon}>m>n^{\delta}$ for some $\epsilon>0$.
Then, for all sufficiently large $n$, the value 
$f(\sigma)>2/n$ can be achieved only on permutations $\sigma$ with 
$p(\sigma)>m$. The number of such permutations $\sigma$ does not exceed
$\displaystyle {n \choose m} (n-m)!=n!/m!$, that is, the probability that 
a random permutation $\sigma$ satisfies $f(\sigma)>2/n$ does not 
exceed $\exp\{-n^{\delta}\}$ for large $n$.

\head 5. The General Case \endhead

It appears that the difference between the general case and the symmetric 
case of Section 4 is not as substantial as the difference between the 
symmetric case and the special cases of Sections 2 and 3.
Our main result is:
\proclaim{(5.1) Theorem} Let 
$f: S_n \longrightarrow {\Bbb R}$ be the function defined by (1.1.1) 
or (1.3.1), let $\overline{f}$ be 
the average value of $f$ on $S_n$ and let $f_0=f-\overline{f}$.
Let $\tau$ be an optimal permutation:
$f_0(\tau)=\max_{\sigma \in S_n} f_0(\sigma)$. Let us choose an 
integer $3 \leq k \leq n-5$ and 
a number $0< \gamma <1$. Let 
$$\beta(n,k)={k-2 \over n^2-nk+k-2}.$$
The probability that a random permutation $\sigma \in S_n$ satisfies 
$$f_0(\sigma) \geq \gamma\beta(k,n) f_0(\tau)$$
is at least
$${(1-\gamma) \beta(k,n) \over 5 k!}.$$
\endproclaim
In particular, by choosing an appropriate $k$, we 
obtain the following corollary.
\proclaim{(5.2) Corollary} 

\roster
\item Let us fix any $\alpha>1$. Then there exists a 
$\delta=\delta(\alpha)>0$ 
such that for all sufficiently large $n \geq N(\alpha)$ the probability that a 
random permutation $\sigma$ in $S_n$ satisfies the inequality
$$f_0(\sigma) \geq {\alpha \over n^2} f_0(\tau)$$
is at least $\delta n^{-2}$. In particular, one can choose 
$\delta=\exp\bigl\{-c\alpha \ln \alpha \bigr\}$ for some absolute 
constant $c>0$.
\item Let us fix any $\epsilon>0$. Then there exists a 
$\delta=\delta(\epsilon)<1$ 
such that for all sufficiently large $n \geq N(\epsilon)$ the probability that 
a random permutation $\sigma$ in $S_n$ satisfies the inequality 
$$f_0(\sigma) \geq n^{-1-\epsilon} f_0(\tau)$$
is at least $\exp\{-n^{\delta}\}$. In particular, one can choose 
any $\delta>1-\epsilon$.
\endroster
\endproclaim
{\hfill \hfill \hfill} \qed

We prove Theorem 5.1 in Section 11. 
As in Section 2, we conclude that for any fixed $\alpha>1$ there is 
a randomized $O(n^2)$ algorithm which produces a permutation $\sigma$ 
satisfying (1). If are willing to settle for an algorithm of 
mildly exponential complexity, we can achieve the bound of type (2), which 
is weaker than the corresponding bound of Corollary 3.2. 

In Section 11, we construct an example of a function of type (1.3.1) with 
an even sharper spike distribution than in example 4.1.

\head 6. Preliminaries \endhead

First, we prove Lemma 1.5.
\demo{Proof of Lemma 1.5} 
Let us choose a pair of indices $1 \leq i \ne j \leq n$. Then, as 
$\sigma$ ranges over the symmetric group $S_n$, the ordered pair 
$\bigl(\sigma(i), \sigma(j)\bigr)$ ranges over all ordered pairs $(k,l)$ 
with $1 \leq k \ne l \leq n$ and each such a pair $(k,l)$ appears $(n-2)!$ 
times. Similarly, for each index $1 \leq i \leq n$, the index $\sigma(i)$ 
ranges over the set $\{1, \ldots, n\}$ and each $j \in \{1, \ldots, n\}$ 
appears $(n-1)!$ times. Therefore, 
$$\split \overline{f}&={1 \over n!} \sum_{\sigma \in S_n} 
\sum_{i,j=1}^n b_{\sigma(i) \sigma(j)} 
a_{ij}=\sum_{i,j=1}^n \Bigl( a_{ij}{1 \over n!}\sum_{\sigma \in S_n} 
b_{\sigma(i) \sigma(j)}\Bigr)\\&=
{1 \over n(n-1)}\sum_{i \ne j}^n a_{ij} \beta_1 + 
{1 \over n}\sum_{i=1}^n a_{ii} \beta_2 =
{\alpha_1 \beta_1 \over n(n-1)} + {\alpha_2 \beta_2 \over n} \endsplit$$
and the proof follows.  
{\hfill \hfill \hfill} \qed
\enddemo
Suppose that 
$f(\sigma)=\langle B, \sigma(A) \rangle$ for some matrices $A$ and $B$ 
and all $\sigma \in S_n$
and suppose that the maximum value of $f$ is attained at a permutation 
$\tau$. Let $A_1=\tau(A)$ and let $f_1(\sigma)=\langle B, \sigma(A_1) 
\rangle$. Then $f_1(\sigma)=f(\sigma \tau)$, hence the maximum value 
of $f_1$ is attained at the identity permutation $e$ and the distribution
of values of $f$ and $f_1$ is the same. We observe that if $A$ is symmetric
then $A_1$ is also symmetric, and if $A$ has constant row and 
column sums and a constant diagonal then so does $A_1$ (see also Section 7).
Hence, as long as the distribution of values of $f$ is concerned, 
without loss of generality we may assume that the maximum of $f$ is attained
at the identity permutation $e$.  
\definition{(6.1) Definition} Let $f: S_n \longrightarrow {\Bbb R}$ be 
a function. Let us define function $g: S_n \longrightarrow {\Bbb R}$ by 
$$g(\sigma)={1 \over n!} \sum_{\omega \in S_n} f(\omega^{-1} \sigma \omega).$$
We call $g$ the {\it central projection} of $f$.
\enddefinition
The following simple observation is quite important for our approach.
\proclaim{(6.2) Lemma} let $f: S_n \longrightarrow {\Bbb R}$ be a 
function such that $f(e) \geq f(\sigma)$ for all $\sigma \in S_n$ and 
let $g$ be the central projection of $f$.
Then $g(e)=f(e) \geq g(\sigma)$ for all $\sigma \in S_n$ and
the average values of $f$ and $g$ are equal: $\overline{f}=\overline{g}$.
\endproclaim
\demo{Proof} We observe that 
$\omega^{-1} e \omega =e$
for all $\omega \in S_n$ and hence $g(e)=f(e)$. 
Moreover, for any $\sigma \in S_n$
$$g(\sigma)={1 \over n!} \sum_{\omega \in S_n} f(\omega^{-1} \sigma \omega) 
\leq {1 \over n!} \sum_{\omega \in S_n} f(e) =g(e).$$
Finally,
$$\split \overline{g} &={1 \over n!} \sum_{\sigma \in S_n} g(\sigma) =
{1 \over n!} \sum_{\sigma \in S_n} {1 \over n!}
\sum_{\omega \in S_n} f(\omega^{-1} \sigma
\omega)={1 \over n!} \sum_{\omega \in S_n} 
\biggl({1 \over n!} \sum_{\sigma \in S_n}
f(\omega^{-1} \sigma \omega) \biggr)\\ &= 
{1 \over n!} \sum_{\sigma \in S_n} \overline{f}=
\overline{f} \endsplit$$ 
and the proof follows.
{\hfill \hfill \hfill} \qed
\enddemo
Moreover, one can observe that the averages of $f$ and $g$ on the 
$k$-th ring $U(e, k)$ coincide for all $k=0, \ldots, n$, see 
Definition 1.2.

We will rely on a Markov type estimate, which asserts, roughly, that 
a function with a sufficiently large average takes sufficiently 
large values sufficiently often.
\proclaim{(6.3) Lemma} Let $X$ be a finite set and let 
$f: X \longrightarrow {\Bbb R}$ be a function. Suppose that $f(x) \leq 1$ 
for all $x \in X$ and that 
$${1 \over |X|} \sum_{x \in X} f(x) \geq \beta \quad 
\text{for some} \quad \beta>0.$$
Then for any $0 < \gamma < 1$ we have
$$
\big|\bigl\{x \in X: f(x) \geq \beta \gamma\bigr\} \big|
\geq \beta (1-\gamma) |X|.$$
\endproclaim 
\demo{Proof} We have 
$$\split &\beta  \leq  {1 \over |X|}
\sum_{x \in X} f(x) ={1 \over |X|}\sum_{x: f(x)< \beta \gamma}
f(x) + {1 \over |X|}
\sum_{x: f(x) \geq \beta \gamma} f(x) \\ & \leq  \beta \gamma
+ {\big|\bigl\{x: f(x) \geq \beta \gamma \bigr\} \big| \over |X|}. \endsplit$$
Hence
$$\big|\{x: f(x) \geq \beta \gamma\}\big| \geq \beta (1-\gamma) |X|.$$
{\hfill \hfill \hfill} \qed
\enddemo
Finally, we need some facts about the structure of the symmetric 
group $S_n$ (see, for example, [6]).
\specialhead (6.4) The conjugacy classes of $S_n$ \endspecialhead 
Let us 
fix a permutation $\rho \in S_n$. As $\omega$ ranges over the symmetric 
group $S_n$, the permutation $\omega^{-1} \rho \omega$ ranges over
the conjugacy class of $X(\rho)$ of $\rho$, that is the set of  
permutations that have the same cycle structure as $\rho$. 

We will be using the following facts.

\subhead{(6.4.1) Central projections and conjugacy classes}
\endsubhead 
If
$f: S_n \longrightarrow {\Bbb R}$ is a function and 
$g: S_n \longrightarrow {\Bbb R}$ its central projection, then 
$$g(\rho)={1 \over |X(\rho)|} \sum_{\sigma \in X(\rho)} 
f(\sigma).$$ 
If $X \subset S_n$ is a set which splits into a union of 
conjugacy classes $X(\rho_i): i \in I$, and for each such a class
we have 
$${1 \over |X(\rho_i)|} \sum_{\sigma \in X(\rho_i)} f(\sigma) \geq \alpha$$ 
for some number $\alpha$, then 
$${1 \over |X|} \sum_{\sigma \in X} f(\sigma) \geq \alpha.$$  
\subhead (6.4.2) Permutations with no fixed points and 2-cycles 
\endsubhead
Let us fix some positive integers $c_i: i=1, \ldots, m$ and let 
$a_n$ be the number of permutations in $S_n$ that have no cycles of 
length $c_i$ for $1 \leq i \leq m$. 
The exponential generating function for $a_n$ is given by
$$\sum_{n=0}^{\infty} {a_n \over n!} x^n ={1 \over 1-x}
 \exp\Bigl\{- \sum_{i=1}^m {x^{c_i} \over c_i} \Bigr\},$$
where we agree that $a_0=1$, see, for example, pp. 170--173 of [7]. 
It follows that the number of permutations $\sigma \in S_n$ without 
fixed points is asymptotically $e^{-1} n!$ and without 
fixed points and 2-cycles is $e^{-3/2} n!$. We will use that the 
first number exceeds $n!/3$ and the second number exceeds $n!/5$ for 
$n \geq 5$.
\subhead (6.4.3) Permutations with many fixed points and 2-cycles
\endsubhead
The number of permutations $\sigma \in S_n$ with at least $k$ 
fixed points is at most $n!/k!$, since to choose such a permutation, 
we can first choose $k$ fixed points in $\displaystyle {n \choose k}$ ways
and then 
choose an arbitrary permutation of the remaining $(n-k)$ elements in 
$(n-k)!$ ways (some permutations will be counted several times). 
Similarly, the number of permutations $\sigma \in S_n$ with at least 
$k$ transpositions (2-cycles) is at most $\displaystyle {n! \over k! 2^k}$,
 since to choose such a permutation,
we first choose some $k$ pairs in 
$\displaystyle {n! \over (n-2k)! k! 2^k}$ ways 
and then an arbitrary permutation of the remaining $n-2k$ elements 
in $(n-2k)!$ ways (again, some permutations will be counted several times).
 
\head 7. Action of the Symmetric Group in the Space of Matrices \endhead

The crucial observation for our approach is that the vector space of all 
central projections $g$ of functions $f$ defined by (1.1.1) or (1.3.1) is 
4-, 3-, or 2- dimensional depending on whether we consider the general case, 
the  
cases of Sections 3 and 4 or the special case of Section 2. If we require, 
additionally, that $\overline{f}=0$ then the dimensions drop by 1 
to 3, 2 and 1, respectively. This fact is explained by the representation 
theory of the symmetric group (see, for example, [6]).
In this section, we review some facts that we need.
Our notation is inspired by the generally accepted notation of the 
representation theory. 

We describe some important invariant subspaces of the action of $S_n$ in 
the space of $n \times n$ matrices $\Mat_n$ by simultaneous permutations 
of rows and columns. We recall that $n \geq 4$.
\subhead(7.1) Subspace $L_n$ \endsubhead Let $L^1_n$ be the space of 
constant matrices $A$:
$$a_{ij}=\alpha\quad \text{for some} \quad \alpha \quad 
\text{and all} \quad 1 \leq i,j  \leq n.$$
Let $L^2_n$ be the subspace of scalar matrices $A$:
$$a_{ij}=\cases \alpha &\text{if \ } i=j \\ 0 &\text{if \ } i \ne j 
\endcases \quad \text{for some} \quad \alpha.$$ 
Finally, Let $L_n=L_n^1 + L_n^2$. One can observe that 
$\dim L_n=2$ and that $L_n$ is the subspace of all matrices that remain 
fixed under the action of $S_n$.
\subhead (7.2) Subspace $L_{n-1,1}$ \endsubhead Let $L_{n-1,1}^1$ be 
the subspace of matrices with identical rows and such that the sum of 
entries in each row is 0:
$$A=\left( \matrix \alpha_1 & \alpha_2 & \ldots & \alpha_n  \\ 
\alpha_1 & \alpha_2 & \ldots & \alpha_n \\ \ldots & \ldots & \ldots & \ldots
\\ \alpha_1 & \alpha_2 & \ldots & \alpha_n \endmatrix \right), 
\quad \text{where} \quad \alpha_1 + \ldots + \alpha_n=0.$$
Similarly, let $L_{n-1,1}^2$ be the subspace of matrices with identical 
columns and such that the sum of entries in each column is 0:
$$A=\left( \matrix \alpha_1 &\alpha_1 &\ldots &\alpha_1 \\ 
\alpha_2 &\alpha_2 &\ldots &\alpha_2 \\ \ldots & \ldots & \ldots & \ldots
\\ \alpha_n &\alpha_n &\ldots &\alpha_n \endmatrix \right), \quad 
\text{where} \quad \alpha_1 + \ldots +\alpha_n=0.$$ 
Finally, let $L_{n-1,1}^3$ be the subspace of diagonal matrices with the zero 
sum on the diagonal:
$$A=\left( \matrix \alpha_1 & 0 & \ldots & 0 & 0 \\ 0 &\alpha_2 &\ldots &0 & 0
 \\
\ldots & \ldots & \ldots & \ldots \\ 0 & 0 & \ldots & 0 & \alpha_n \endmatrix
\right), \quad \text{where} \quad \alpha_1 + \ldots + \alpha_n=0.$$ 
Let $L_{n-1,1}=L_{n-1,1}^1 + L_{n-1,1}^2+ L_{n-1,1}^3$. One can check that
the dimension of each of $L_{n-1,1}^1$, $L_{n-1,1}^2$ and 
$L_{n-1,1}^3$ is $n-1$ and that 
$\dim L_{n-1,1}=3n-3$. Moreover, the subspaces 
$L_{n-1,1}^1$, $L_{n-1,1}^2$ and $L_{n-1,1}^3$ 
do not contain non-trivial invariant subspaces.
The action of $S_n$ in $L_{n-1,1}$, although 
non-trivial, is not very complicated. One can show that if 
$A \in L_{n-1,1}+L_n$, then the problem of optimizing $f(\sigma)$ defined 
by (1.1.1) reduces to the Linear Assignment Problem.
\subhead (7.3) Subspace $L_{n-2,2}$ \endsubhead Let us define $L_{n-2,2}$ as 
the subspace of all {\it symmetric} matrices $A$ with row and 
column sums equal to 0 and zero diagonal 
$$\split &a_{ij}=a_{ji} \quad \text{for all} \quad 1 \leq i,j \leq n; \\
& \sum_{i=1}^n a_{ij}=0 \quad \text{for all} \quad j=1, \ldots, n;  \\ 
& \sum_{j=1}^n a_{ij}=0 \quad \text{for all} \quad i=1, \ldots, n 
\quad \text{and} \\ 
&a_{ii}=0 \quad \text{for all} \quad i=1, \ldots, n. \endsplit$$
One can check that $L_{n-2,2}$ is an invariant subspace and that 
$\dim L_{n-2,2}=(n^2-3n)/2$.
Besides, $L_{n-2,2}$ contains no non-trivial invariant subspaces.
\subhead (7.4) Subspace $L_{n-2,1,1}$ \endsubhead Let us define $L_{n-2,1,1}$
as
the subset of all {\it skew symmetric} matrices $A$ with row and 
column sums equal to 0:
$$\split &a_{ij}=-a_{ji} \quad \text{for all} \quad 1 \leq i,j \leq n; \\
& \sum_{i=1}^n a_{ij}=0 \quad \text{for all} \quad j=1, \ldots, n \quad 
\text{and}  \\ 
& \sum_{j=1}^n a_{ij}=0 \quad \text{for all} \quad i=1, \ldots, n. \endsplit$$
One can check that $L_{n-2,1,1}$ is an invariant subspace and that 
$\dim L_{n-2,1,1}=(n^2-3n)/2+1$. Similarly, $L_{n-2,1,1}$ contains 
no non-trivial invariant subspaces.
\bigskip
One can check that $\Mat_n=L_n+L_{n-1,1}+L_{n-2,2}+L_{n-2,1,1}$.
The importance of the subspaces (7.1)--(7.4) is explained by the fact that
they are the {\it isotypical components} of the irreducible representations 
of the symmetric group in the space of matrices. The following proposition
follows from the representation theory of the symmetric group [6].
\proclaim{(7.5) Proposition} For an $n \times n$ matrices $A$ and $B$,
where $n \geq 4$, let $f: S_n \longrightarrow {\Bbb R}$ be the function 
defined by (1.1) and let $g: S_n \longrightarrow {\Bbb R}$,
$$ g(\sigma)={1 \over n!} \sum_{\omega \in S_n} 
f\bigl(\omega^{-1} \sigma \omega\bigr)$$
be the central projection of $f$.
Given a permutation $\sigma \in S_n$, let 
$$p(\sigma)=\big|\{i: \sigma(i)=i \}\big| \quad \text{and} \quad 
t(\sigma)=\big|\{ i<j: \sigma(i)=j \ \text{and} \ \sigma(j)=i \}\big|$$
be the number of fixed points of the permutation and the number of 
2-cycles in the permutation correspondingly.
\roster
\item If $A \in L_n$ then $g$ is a scalar multiple of the constant function
$$\chi_n(\sigma)=1 \quad \text{for all} \quad \sigma \in S_n;$$
\item If $A \in L_{n-1,1}$ then $g$ is a scalar multiple of the function
$$\chi_{n-1,1}(\sigma)=p(\sigma)-1
\quad \text{for all} \quad \sigma \in S_n; $$
\item If $A \in L_{n-2,2}$ then $g$ is a scalar multiple of the 
function 
$$\chi_{n-2,2}(\sigma)=t(\sigma)+{1 \over 2} p^2(\sigma)-{3 \over 2} 
p(\sigma) \quad
\text{for all} \quad \sigma \in S_n;$$
\item If $A \in L_{n-2,1,1}$ then $g$ is a scalar multiple of 
the function 
$$\chi_{n-1,1,1}(\sigma) ={1 \over 2}p^2(\sigma)-{3 \over 2}p(\sigma) 
-t(\sigma)+1 \quad
\text{for all} \quad \sigma \in S_n.$$ 
\endroster 
\endproclaim
The functions $\chi_n, \chi_{n-1,1}, \chi_{n-2,2}$ and $\chi_{n-1,1,1}$ are 
the characters of corresponding irreducible representations of 
$S_n$ for $n \geq 4$. They are linearly independent, and, moreover 
orthogonal: $\sum_{\sigma \in S_n} \chi_i(\sigma) \chi_j(\sigma)=0$
for two characters of different irreducible representation of $S_n$.
In particular,
$$\sum_{\sigma \in S_n} \chi_{n-1,1}(\sigma)=
\sum_{\sigma \in S_n} \chi_{n-2,2}(\sigma)=
\sum_{\sigma \in S_n} \chi_{n-2,1,1}(\sigma)=0,$$
hence the average value of all but the trivial character $\chi_n$ is 0.
\remark{(7.6) Remark} It follows [6] that each of the 
functions $\chi_n, \chi_{n-1,1}, \chi_{n-2,2}$ and \break $\chi_{n-2,1,1}$ is
the objective function (1.3.1) in some generalized problem with 
a tensor $c^{ij}_{kl}$ (see Section 1.3)
with the property that for all $k$ and $l$ the matrix $A=(a_{ij})$
for $a_{ij}=c^{ij}_{kl}$
belongs to the corresponding subspace. 
Since the set of all functions (1.3.1) is closed under linear combinations, 
it follows that every function $f \in \spa\{\chi_n, \chi_{n-1,1}, \chi_{n-2,2},
\chi_{n-2,1,1} \}$ is an objective function in the generalized 
problem.
\endremark

\head 8. The Bullseye Case. Proofs \endhead

In this section, we prove Theorem 2.1 and Theorem 2.3.
An important observation is that 
$A$ satisfies the conditions of Section 2 if and only if 
$A \in L_n +L_{n-2,2}$ (see Section 7). 
\demo{Proof of Theorem 2.1} Without loss of generality, we may assume that
the maximum of $f_0(\sigma)$ is attained at the 
identity permutation $e$ (see Section 6). 
Excluding the non-interesting  case of $f_0 \equiv 0$, by scaling $f$, if 
necessary, we can assume that $f_0(e)=1$. Let $g$ be the central projection
of $f_0$. Then by Lemma 6.2,  $\overline{g}=0$ and 
$1=g(e) \geq g(\sigma)$ for all $\sigma \in S_n$.
Moreover, since $A \in L_n + L_{n-2,2}$,
by Parts 1 and 3 of Proposition 7.5, $g$ must be a linear combination of 
the constant function $\chi_n$ and $\chi_{n-2,2}$. Since $\overline{g}=0$, 
$g$ should be proportional to $\chi_{n-2,2}$ and 
since $g(e)=1$, we have 
$$g={2 \over n^2-3n} \chi_{n-2,2}={2 t+p^2-3p \over
n^2-3n}.$$ 
Now $\sigma \in U(e, k)$ if and only if $p(\sigma)=n-k$. Hence 
$g(\sigma) \geq \alpha(n,k)$ for all $\sigma \in U(e,k)$.
The set $U(e,k)$ splits into disjoint union of conjugacy classes $X(\rho)$ 
and, using (6.4.1), we conclude that for each such $X(\rho)$
$$g(\rho)={1 \over |X(\rho)|} \sum_{\sigma \in X(\rho)} f_0(\sigma) 
\geq \alpha(n,k)$$ and, therefore,
$${1 \over |U(n,k)|} \sum_{\sigma \in U(n,k)} f_0(\sigma) \geq \alpha(n,k),$$
hence the proof follows.
{\hfill \hfill \hfill} \qed
\enddemo 
Using estimates of (6.4.2), one can show that the input of the number 
of 2-cycles $t(\sigma)$ into the average of 
$f_0$ over $U(e, k)$ is asymptotically
negligible, so there is an ``almost equality'' in the formula
of Theorem 2.1.

By estimating the cardinality of the $k$-th ring $U(\tau, k)$, we 
deduce Theorem 2.3.
\demo{Proof of Theorem 2.3} As in the proof of Theorem 2.1, we assume  
that the maximum value of $f_0$ is equal to 1.

Let us estimate the cardinality 
$|U(\tau, n-k)|=|U(e, n-k)|$. Since $\sigma \in U(e,n-k)$ if and only if
$\sigma$ has $k$ fixed points, to choose a $\sigma \in U(e,n- k)$ one 
has to choose $k$ points in 
$\displaystyle{n \choose k}$ ways and then choose a 
permutation of the remaining $n-k$ points without fixed points.
Using (6.4.2), we get
$$|U(\tau, n-k)| \geq  {n \choose k} (n-k)!/3 ={n! \over 3 k!}.$$ 
Applying Lemma 6.3 with $\beta=\beta(n,k)$ and $X=U(\tau, n-k)$, from 
Theorem 2.1, we conclude that
$$\split \PP\Bigl\{\sigma \in S_n: f_0(\sigma) \geq \gamma \beta(n,k)\Bigr\}
 &
 \geq 
{(1-\gamma)\beta(n,k) |U(\tau,n-k)| \over n!} \\ 
&\geq {(1-\gamma)\beta(n,k) \over 3 k!}. \endsplit$$
{\hfill \hfill \hfill} \qed
\enddemo

\head 9. The Pure Case. Proofs \endhead

In this case, $A \in L_n + L_{n-2,1,1} +L_{n-2,2}$ (see Section 7).
As in Section 8, the $L_n$ component contributes just a constant to $f$.
Since the $L_{n-1,1}$ component attributed to the Linear Assignment 
Problem (see Section 7.2) is absent, we call this case ``pure''.

We choose a more convenient basis $g_1$ and $g_2$ in the vector space 
spanned by $\chi_{n-2,2}$ and $\chi_{n-2,1,1}$, namely:
$$g_1=\chi_{n-2,2} + \chi_{n-2,1,1}=p^2-3p+1 \quad 
\text{and} \quad g_2=\chi_{n-2,1,1}-\chi_{n-2,2}=1-2t.$$
\definition{(9.1) Definition} Let $K_p$ (where $p$ stands for 
``pure'') be the set of all functions 
$g: S_n \longrightarrow {\Bbb R}$ such that $g \in \spa\{g_1, g_2\}$,
where $g_1=p^2-3p+1$ and $g_2=1-2t$ and $g(e) \geq g(\sigma)$ 
for all $\sigma \in S_n$, where $e$ is the identity permutation.
We call $K_p$ the {\it central cone}.
\enddefinition
Identifying $\spa\{g_1, g_2\}$ with two-dimensional vector space ${\Bbb R}^2$
(plane), we see that the conditions $g(e) \geq g(\sigma)$ define the 
central cone $K$ as 
a convex cone in ${\Bbb R}^2$. Our goal is to find the extreme 
rays $r_1$ and $r_2$ of $K$, so that every function $g \in K$ can be 
written as a non-negative linear combination of $r_1$ and $r_2$.

First, we prove a useful technical result.

\proclaim{(9.2) Lemma} For a permutation $\sigma \in S_n$, $\sigma \ne e$,
let $a_{\sigma} \in {\Bbb R}^2$ be the point
$$a_{\sigma}=\Bigl(p(\sigma), \quad {2t(\sigma) \over n-p(\sigma)} \Bigr).$$
Let $P=\conv\bigl\{a_{\sigma}: \sigma \ne e\bigr\}$ be 
the convex hull 
of all such points $a_{\sigma}$. 

If $n$ is even, the extreme points of $P$ are
$$(0,0), \quad (n-3, 0), \quad (n-2, 1) \quad \text{and} \quad  (0,1).$$

If $n$ is odd, the extreme points of $P$ are 
$$(0,0), \quad (n-3, 0), \quad (n-2,1), \quad \bigl(0, (n-3)/n\bigr)
\quad \text{and} \quad \quad (1, 1).$$
\endproclaim
\demo{Proof} The set of all possible values $\bigl(p(\sigma), t(\sigma)\bigr)$,
where $\sigma \ne e$,
consists of all pairs of non-negative integers $(p,t)$ such that 
$p \leq n-2$, $2t \leq n$ and, additionally, $p+2t \leq n-3$ or 
$p+2t=n$. To find the extreme points of the set of 
feasible points $\bigl(p, 2t/(n-p)\bigr)$, we 
choose a generic vector $(\gamma_1, \gamma_2)$ and 
investigate for which values of $p$ and $t$ the maximum of 
$$\gamma_1 p + \gamma_2 {2 t \over n-p}$$
is attained.

Clearly, we can assume that $\gamma_2 \ne 0$. If $\gamma_2<0$ then
we should choose the smallest possible $t$ which would be $t=0$ unless
$p=n-2$ when we have to choose $t=1$. 
Depending on the sign of $\gamma_1$, this produces the 
following pairs 
$$(p,t)=\Bigl\{(0, 0),\quad  (n-3, 0), \quad (n-2, 1) \Bigr\}.$$  
If $\gamma_2>0$ then the largest possible value of $2t/(n-p)$ is 1.
If $\gamma_1>0$ this produces the (already included) point
$$(p,t)=(n-2, 1\bigr).$$
If $\gamma_1<0$ we get 
$$(p,t)=(0,n/2) \quad \text{for even\ } n$$
and
$$ (p,t)=\Bigl\{\bigl(0, (n-3)/2 \bigr), \bigl(1, (n-1)/2\bigr)
\Bigr\} \quad \text{for odd} \quad n.$$
Summarizing, the extreme points of $P$ are 
$$(0,0), \quad (n-3, 0), \quad (n-2, 1), \quad (0,1) \qquad 
\text{for even\ } n$$
and
$$(0,0), \quad (n-3, 0), \quad (n-2,1), \quad \bigl(0, (n-3)/n\bigr), 
\quad (1, 1) \qquad \text{for odd\ }n$$
as claimed.
{\hfill \hfill \hfill} \qed
\enddemo
Now we describe the central cone $K_p$.
\proclaim{(9.3) Lemma}
For $n \geq 4$ let us define the functions $r_1, r_{2e}$ and $r_{2o}:
S_n \longrightarrow {\Bbb R}$ by 
$$\split &r_1=1-2t, \\ 
         &r_{2e}={p^2-3p-n-6t+2tn+4 \over n^2-4n+4} \quad 
\text{and} \\ 
&r_{2o}={p^2-3p-n-4t+2tn+3 \over n^2-4n+3}.
\endsplit$$  
Then 
\roster
\item If $n$ is even then
 $K_p$ is a 2-dimensional convex cone with the extreme rays 
spanned by $r_1$ and $r_{2e}$;
\item If $n$ is odd then 
$K_p$ is a 2-dimensional convex cone with the extreme rays 
spanned by $r_1$ and $r_{2o}$. Cone $K_p$ contains $r_{2e}$;
\item If $e \in S_n$ is the identity, then
$$r_1(e)=r_{2e}(e)=r_{2o}(e)=1.$$
\endroster
\endproclaim
\demo{Proof} A function $g \in K_p$ can be written as a linear combination
$g=\alpha_1 g_1 +\alpha_2 g_2$. Since $p(e)=n$ and $t(e)=0$, we have 
$g(e)=\alpha_1(n^2-3n+1) +\alpha_2$. Therefore, the inequalities 
$g(e) \geq g(\sigma)$ can be written as 
$$\alpha_1(n^2-3n+1) +\alpha_2 \geq 
\alpha_1\bigl(p(\sigma)^2-3p(\sigma)+1\bigr) +
\alpha_2\bigl(1-2t(\sigma)\bigr),$$
which, for $g \ne e$ is equivalent to 
$$\alpha_1\bigl(n+p(\sigma)-3\bigr) + \alpha_2{ 2t(\sigma) \over n-p(\sigma)} 
\geq 0.$$
Using Lemma 9.2, we conclude that for even $n$, the system is equivalent 
to 
$$\split &\alpha_1 \geq 0 \\
         &(n-3) \alpha_1 + \alpha_2 \geq 0 \endsplit \tag9.3.1$$
and for odd $n$, the system is equivalent to 
$$\split &\alpha_1 \geq 0 \\ 
         &(n-2)\alpha_1 + \alpha_2 \geq 0. \endsplit \tag9.3.2$$
Consequently, every solution $(\alpha_1, \alpha_2)$ of (9.3.1) is 
a non-negative linear combination of $(0,1)$ and $(1, 3-n)$ and 
every solution of (9.3.2) is a non-negative linear combination of 
$(0,1)$ and $(1,2-n)$.
  
The functions $r_1, r_{2e}$ and $r_{2o}$ are obtained from 
$g_2, g_1+(3-n)g_2$ and $g_1 +(2-n)g_2$ respectively by scaling 
so that the value at the identity becomes equal to 1.

Since every solution of (9.3.1) is a solution of (9.3.2), we conclude 
that $r_{2e} \in K_p$ for odd $n$ as well.
{\hfill \hfill \hfill} \qed
\enddemo
\remark{(9.4) Remark} If $n$ is even, then 
$r_{2o} \notin K_p$. Indeed, if $\sigma$ is a product of $n/2$ 
commuting transpositions, so that $p(\sigma)=0$ and $t(\sigma)=n/2$, then 
$r_{2o}(\sigma)=(n^2-3n+3)/(n^2-4n+3)>1=r_{2o}(e)$. 
$$\epsffile{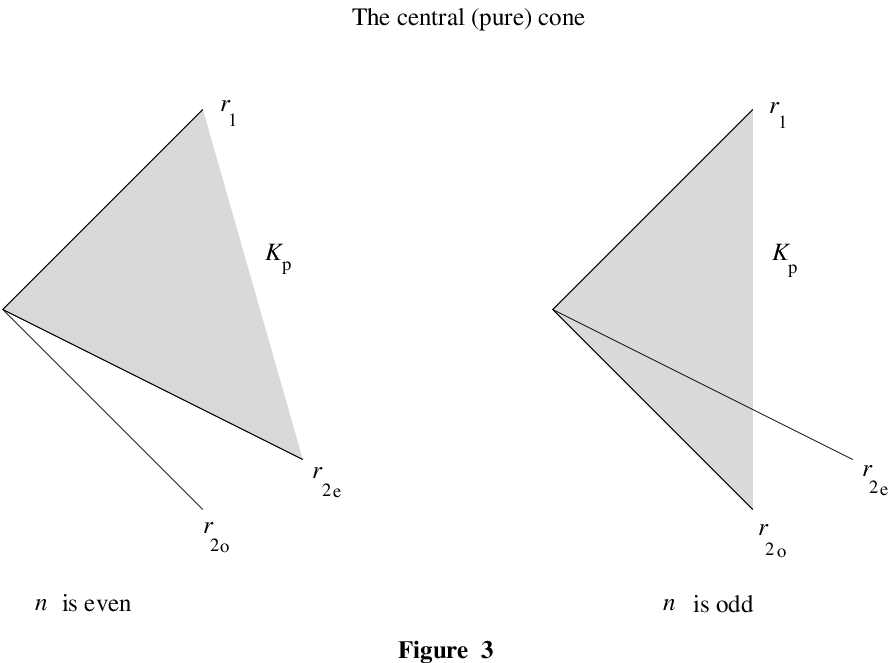}$$
The functions $r_{2o}$ and $r_{2e}$ have the bullseye distribution of 
Section 2. The distribution type of $r_1$ may be characterized as that of 
a ``damped oscillator'' with the averages over the $k$-ring 
$U(e,k)$ changing sign and going fast to 0 as $k$ grows.
Hence a typical function from the central cone has a ``weak'' bullseye 
type distribution, which becomes weaker as the function becomes closer 
to $r_1$.
\endremark
\proclaim{(9.5) Lemma} Let $g \in K_p$ be a function such that 
$g(e)=1$. For any $1 \leq k \leq n-2$, let $\sigma_k$ be a permutation 
such that $p(\sigma_k)=k$ and $t(\sigma_k)=0$ and let $\theta_k$ be 
a permutation such that $p(\theta_k)=k$ and $t(\theta_k)=1$. 
Then 
$$\max\bigl\{g(\sigma_k), g(\theta_k) \bigr\} \geq {k^2-3k+1 \over n^2-3n+1}.$$
\endproclaim
\demo{Proof} Applying Lemma 9.3, we may assume that 
$g$ is a convex combination of $r_1$ and $r_{2o}$, hence
$g=\alpha_1 r_1 +\alpha_2 r_{2o}$,
for some $\alpha_1, \alpha_2 \geq 0$ and $\alpha_1 + \alpha_2=1$.
Then 
$$g(\sigma_k)=\alpha_1 + \alpha_2 {k^2-3k -n +3 \over n^2-4n+3}$$  
and
$$g(\theta_k)=-\alpha_1 + \alpha_2 
{k^2-3k+n-1 \over n^2-4n+3}.$$  
We observe that if $\alpha_1=1$ and $\alpha_2=0$ then 
$g(\sigma_k)>g(\theta_k)$ and if $\alpha_1=0$ and $\alpha_2=1$ 
then $g(\sigma_k)<g(\theta_k)$. Moreover, as $(\alpha_1, \alpha_2)$ 
change from $(1,0)$ to $(0,1)$ function $g(\sigma_k)$ decreases and
function $g(\theta_k)$ increases. Hence the minimum of 
$\max\{g(\sigma_k), g(\theta_k) \}$ is attained when 
$g(\sigma_k)=g(\theta_k)$. This produces the system of linear 
equations
$$\alpha_1 + \alpha_2 {k^2-3k -n +3 \over n^2-4n+3}=-\alpha_1 + \alpha_2 
{k^2-3k+n-1 \over n^2-4n+3}$$
and $$\alpha_1 +\alpha_2=1$$
with the solution
$$\alpha_1={n-2 \over n^2-3n+1} \quad \text{and} \quad 
\alpha_2={n^2-4n+3 \over n^2-3n+1}.$$
The corresponding value of $g(\theta_k)=g(\sigma_k)$ is 
$${k^2-3k+1 \over n^2-3n+1},$$
which completes the proof.
{\hfill \hfill \hfill} \qed 
\enddemo
Now we are ready to prove Theorem 3.1.
\demo{Proof of Theorem 3.1} Without loss of generality, we may 
assume that the maximum value of $f_0$ is attained at the identity 
permutation $e$ (see Section 6). Excluding an obvious case of 
$f_0 \equiv 0$, by scaling $f$, if necessary, we may assume 
that $f_0(e)=1$. Let $g$ be the central projection of $f_0$. By Lemma 6.2, 
$g(e)=f_0(e)=1 \geq g(\sigma)$ for all $\sigma \in S_n$ and 
$\overline{g}=0$. Moreover, since $A \in L_n + L_{n-2,2} + L_{n-2,1,1}$,
by Proposition 7.5, $g$ must be a linear combination of the constant 
function $\chi_n$ and functions $\chi_{n-2,2}$ and $\chi_{n-2,1,1}$.
Since $\overline{g}=0$, $g$ is a linear combination of 
$\chi_{n-2,2}$ and $\chi_{n-2,1,1}$ alone. Therefore, $g$ lies in the 
central cone: $g \in K_p$, see Definition 9.1.

Let us choose a $3 \leq k \leq n-5$ and let $X_k$ be the set 
of permutations $\sigma$ such that $p(\sigma)=k$ and $t(\sigma)=0$ and 
let $Y_k$ be the set of permutations $\theta$ such that 
$p(\theta)=k$ and $t(\theta)=1$. To choose a permutation 
$\sigma \in X_k$, one has to choose $k$ fixed points in 
$\displaystyle {n \choose k}$ ways and then a permutation without fixed 
points or 2-cycles on the remaining $(n-k)$ points. Then, by (6.4.2)
$$|X_k| \geq {1 \over 5} {n \choose k} (n-k)! =
 {1 \over 5} {n! \over k!}.$$  
Similarly, to choose a permutation $\theta \in Y_k$, one has to choose 
a 2-cycle in $\displaystyle {n \choose 2}$ ways, $k$ fixed points 
in $\displaystyle {n-2 \choose k}$ ways and a permutation without 
fixed points or 2-cycles on the remaining $(n-k-2)$ points. Then, by 
(6.4.2)
$$|Y_k| \geq {1 \over 5}{n \choose 2} {n-2 \choose k} (n-k-2)!=
{n! \over 10k!}.$$ 
Let us choose a permutation $\sigma \in X_k$ and a permutation 
$\theta \in Y_k$ and let $Z=X_k$ if $g(\sigma_k) \geq g(\theta_k)$ and 
$Z=Y_k$ otherwise. Then 
$$|Z| \geq  {n! \over 10k!}$$
and by Lemma 9.5,
$$g(\sigma) \geq {k^2-3k+1 \over n^2-3n+1} \quad \text{for all} 
\quad \sigma \in Z.$$
The set $Z$ is a disjoint union of some conjugacy classes $X(\rho)$ 
and for each $X(\rho)$ by (6.4.1), we have 
$$g(\rho)={1 \over |X(\rho)|} \sum_{\sigma \in X(\rho)} f_0(\sigma) \geq 
{k^2-3k+1 \over n^2-3n+1}$$ 
and hence 
$${1 \over |Z|} \sum_{\sigma \in X(\rho)} f_0(\sigma) \geq 
{k^2-3k+1 \over n^2-3n+1}.$$ 
Applying Lemma 6.3 with $X=Z$ and $\beta=\beta(n,k)$, we 
get that 
$$\PP\Bigl\{\sigma \in S_n: f_0(\sigma) \geq \gamma
\beta(n,k) \Bigr\} \geq {(1-\gamma)\beta(n,k) \over 10k!}.$$ 
\enddemo

\head 10. The Symmetric Case. Proofs \endhead

In this case, $A \in L_n + L_{n-1,1} + L_{n-2,2}$ (see Section 7). As 
in Sections 8 and 9, the $L_n$ component contributes a just a constant to 
$f$.   
We choose a more convenient basis $g_1$ and $g_2$ in the vector 
space spanned by $\chi_{n-1,1}$ and $\chi_{n-2,2}$, namely
$$g_1=\chi_{n-1,1}=p-1 \quad \text{and} \quad g_2=2\chi_{n-2,2}+
3\chi_{n-1,1}=p^2+2t-3,$$
where $p(\sigma)$ is the number of fixed points of 
$\sigma$ and $t(\sigma)$ is the number of 2-cycles in $\sigma$.
\definition{(10.1) Definition} Let $K_s$ ( where $s$ stands 
for ``symmetric'') be the set of all functions 
$g: S_n \longrightarrow {\Bbb R}$ such that $g \in \spa\{g_1, g_2\}$,
where $g_1=p-1$ and $g_2=p^2+2t-3$ and $g(e) \geq g(\sigma)$ 
for all $\sigma \in S_n$, where $e$ is the identity permutation.
We call $K$ the {\it central cone}.
\enddefinition
Identifying $\spa\{g_1, g_2\}$ with two-dimensional vector space ${\Bbb R}^2$
(plane), we see that the conditions $g(e) \geq g(\sigma)$ define the 
central cone $K_s$ as 
a convex cone in ${\Bbb R}^2$. Our immediate goal is to find the extreme 
rays $r_1$ and $r_2$ of $K_s$, so that every function $g \in K_s$ can be 
written as a non-negative linear combination of $r_1$ and $r_2$.
\proclaim{(10.2) Lemma}
For $n \geq 4$ let us define the functions $r_1, r_{2e}$ and $r_{2o}:
S_n \longrightarrow {\Bbb R}$ by 
$$\split &r_1={2np-2n-p^2-3p-2t+6 \over n^2-5n+6}, \\ 
         &r_{2e}={-np+n+p^2+p+2t-4 \over 2n-4} \quad 
\text{and} \\ 
&r_{2o}={-n^2p+np^2+n^2+np+2nt -4n -3p+3 \over 2n^2-7n+3}.
\endsplit$$  
Then 
\roster
\item If $n$ is even then
 $K_s$ is a 2-dimensional convex cone with the extreme rays 
spanned by $r_1$ and $r_{2e}$;
\item If $n$ is odd then 
$K_s$ is a 2-dimensional convex cone with the extreme rays 
spanned by $r_1$ and $r_{2o}$. Cone $K_s$ contains $r_{2e}$;
\item If $e \in S_n$ is the identity, then
$$r_1(e)=r_{2e}(e)=r_{2o}(e)=1.$$
\endroster
\endproclaim
\demo{Proof} A function $g \in K_s$ can be written as a linear combination
$g=\alpha_1 g_1 + \alpha_2 g_2$. Since $p(e)=n$ and $t(e)=0$, we have 
$g(e)=\alpha_1 (n-1) + \alpha_2 (n^2-3)$. 
Therefore, the inequalities $g(e) \geq g(\sigma)$ can be written 
as 
$$\alpha_1 (n-1) +\alpha_2(n^2-3) \geq 
\alpha_1 \bigl(p(\sigma)-1\bigr) +\alpha_2 \bigl(p^2(\sigma)
+2t(\sigma)-3\bigr),$$
which, for $\sigma \ne e$, is equivalent to
$$\alpha_1 +\alpha_2\Bigl(n+p(\sigma)-{2t(\sigma) \over n-p(\sigma)}
\Bigr) \geq 0. \tag10.2.1$$
Applying Lemma 9.2, we observe that (10.2.1) 
is equivalent to the system of two inequalities:
$$\alpha_1 + (2n-3) \alpha_2 \geq 0$$
and
$$\split &\alpha_1 + (n-1) \alpha_2 \geq 0 \quad \text{if \ } n \text{\ is 
even}, \\
&n \alpha_1 +(n^2-n+3) \alpha_2 \geq 0 
 \quad \text{if \ } n \text{\ is odd}. \endsplit$$
Thus every pair $(\alpha_1, \alpha_2)$ satisfying (10.2.1) can be 
written as a non-negative linear combination of 
$(2n-3, -1)$ and $(1-n, 1)$ when $n$ is even and $(2n-3,-1)$ and 
$(-n^2+n-3, n)$ when $n$ is odd.

The generators $r_1, r_{2e}$ and $r_{2o}$ are obtained from 
$(2n-3)g_1-g_2$, $(1-n)g_1 +g_2$ and $(-n^2+n-3)g_1 +ng_2$ 
respectively by scaling so that the value at the identity becomes 
equal to 1. 

It remains to check that $r_{2e} \in K$ for $n$ odd as well. Indeed,
using that $2t+p \leq n$ we have 
$$\split (2n-4)\bigl(r_{2e}-1\bigr) &=
-n(p-1)+p(p+1)+2t-4 -2n+4\\ &=-n(p+1)+p(p+1)+2t 
(p+1)(-n+p)+2t \\ &\leq (p+1)(-n+p)+n-p =p(-n+p) \leq 0. \endsplit$$
{\hfill \hfill \hfill} 
\qed
\enddemo
\remark{(10.3) Remark}
The average value of $r_1$, $r_{2e}$ and $r_{2o}$ on $S_n$ is 0.

The function $r_1: S_n \longrightarrow {\Bbb R}$ 
provides an example of the ``bullseye'' distribution (see Section 2.2).
The maximum value of 1 is attained at the identity and at any transposition.
The positive values of $r_1$ occur on permutations with at least two 
fixed points and $r_1(\sigma)=\Omega\bigl(p(\sigma)/n\bigr)$ if 
$p(\sigma) \geq 3$. 

In contrast, $r_{2e}$ and $r_{2o}$ exhibit a spike type distribution of
Section 4.1. The maximum value of 1 is attained 
at the identity and, for $r_{2e}$, on the product of $n/2$ transpositions, 
or, for $r_{2o}$, on the product of $(n-3)/2$ transpositions.
On the other hand, no 
permutation other than $e$ with at least 2 fixed points yields a positive 
value. 

One can observe that if $n$ is even then $r_{2o} \notin K$. Indeed, if 
$\sigma$ is a product of $n/2$ transpositions then 
$r_{2o}(\sigma)=(2n^2-4n+3)/(2n^2-7n+3) > 1$. 
$$\epsffile{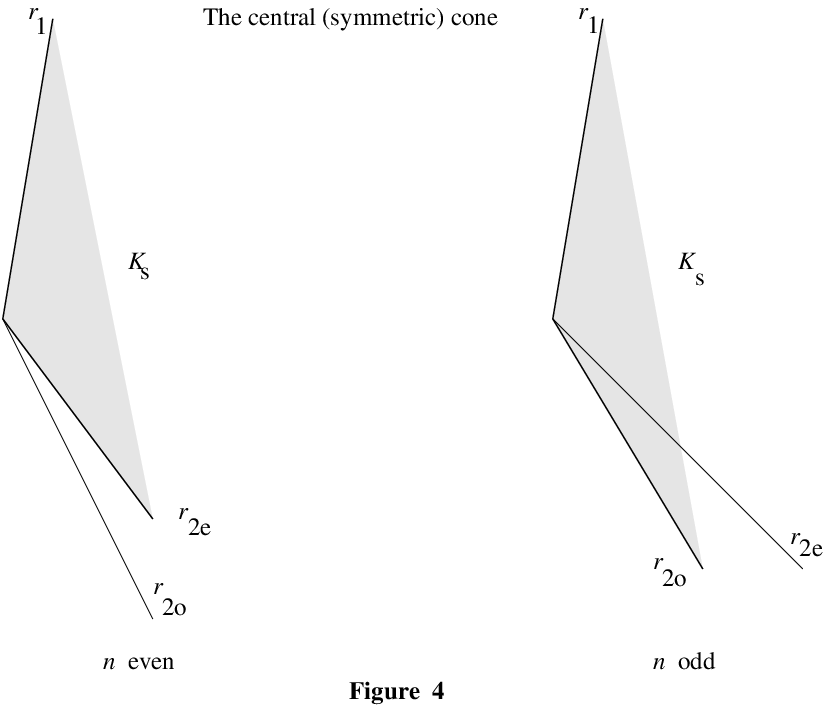}$$
The picture of $K_s$ is very similar to that of $K_p$, see  
Section 9.4.
\endremark
\remark{(10.4) Remark. The spike distribution} Let us consider Example 4.1.
It is seen that 
$f(\sigma)=2b_{\sigma(1) \sigma(2)}$ and hence the maximum value of 
$f$ is indeed 1 and obtained, in particular, on the identity permutation $e$.
Applying Lemma 1.5, we get
$$\overline{f}={1 + 4\gamma(n-2) +0.5(n-2)(n-3) \over n(n-1)}=0.$$
Let us prove that the central projection of $f$ is the function $r_{2e}$ 
of Lemma 10.2. 

Suppose that $g$ is the central projection of $f$. It follows 
that $g$ can be written as a linear combination 
$g=\alpha_1 r_1 + \alpha_2 r_{2e}$.
Since $g(e)=r_1(e)=r_{2e}(e)=1$, we must have $\alpha_1 + \alpha_2=1$.
 Let $\theta=(12)$ be a transposition,
hence $p(\theta)=n-2$ and $t(\theta)=1$. Then $r_1(\theta)=1$ and 
$r_{2e}(\theta)=0$, hence $\alpha_1=g(\theta)$. 

Denoting by $X$ the set of all transpositions in $S_n$, by (6.4.1) we 
get 
$$\split g(\theta) &={1 \over |X|}\sum_{\sigma \in X} f(\sigma)=
2{n \choose 2}^{-1} \sum_{\sigma \in X} b_{\sigma(1) \sigma(2)}\\
&= 
{n \choose 2}^{-1} \Bigl(1+ {n-2 \choose 2} + 4(n-2)\gamma \Bigr)=0.
\endsplit$$
Therefore, $g=r_{2e}$. Let $X(\rho)$ be a conjugacy class with 
$p(\rho)=n-k$. Then 
$$\split{1 \over |X(\rho)|} \sum_{\sigma \in X(\rho)} f(\sigma)
=r_{2e}(\rho) &={-n(n-k)+n+(n-k)^2+(n-k)+2t(\rho)-4 \over 2n-4} \\ 
  &\leq {-nk+k^2-k+3n-4 \over 2n-4}. 
\endsplit$$ 
Since the $k$-th ring $U(e,k)$ splits into a disjoint union of 
conjugacy classes $X(\rho)$ with $p(\rho)=n-k$, we conclude by (6.4.1) 
that 
$${1 \over |U(e,k)|} \sum_{\sigma \in U(e,k)} f(\sigma) \leq 
{-nk+k^2-k+3n-4 \over 2n-4}$$
as claimed.

More generally, one can prove that for any function $g \in K_s$ 
there is a function $f$ of type (1.1.1) with symmetric $A$, such that
$\overline{f}=0$, $f$ attains its maximum at the identity and 
the central projection of $f$ is $g$.
\endremark 
\remark{(10.5) Remark. Scarcity of relatively good values}
Let us consider the function $f$ of Example 4.2.
We observe that 
$$f=\alpha_1 r_1 + \alpha_2 r_{2e}$$
for
$$\alpha_1={n^2-nm-4n+3m+3 \over n^2-nm+m-3} \quad \text{and} 
\quad \alpha_2={4n-2m-6 \over n^2-nm+m-3}.$$
Thus $f$ is a convex combination of $r_1$ and $r_{2e}$, hence 
$1=f(e) \geq f(\sigma)$ for all $\sigma \in S_n$ and $\overline{f}=0$.
Remark 7.6 implies that $f$ is a generalized function (1.3.1) of the 
required type.
\endremark

\head 11. The General Case. Proofs \endhead

In this Section, we prove Theorem 5.1 and describe 
the ``extreme'' distributions. 

Let us choose a convenient basis in $\spa\{\chi_{n-1,1}, \chi_{n-2,2},
\chi_{n-2,1,1}\}$:
$$g_1=\chi_{n-1,1}=p-1, \quad g_2=\chi_{n-2,2} +\chi_{n-2,1,1}+3\chi_{n-1,1}=
p^2-2 \quad \text{and} $$
$$g_3=\chi_{n-2,1,1}-\chi_{n-2,2}=1-2t.$$
\definition{(11.1) Definition} Let $K$ be the set of all functions
$g \in \spa\{g_1,g_2,g_3\}$ such that 
$g(e) \geq g(\sigma)$ for all $\sigma \in S_n$. We call $K$ the 
{\it central cone}.
\enddefinition
Identifying $\spa\{g_1, g_2, g_3\}$ with a 3-dimensional vector space 
${\Bbb R}^3$, we see that conditions $g(e) \geq g(\sigma)$ define the central 
cone $K$ as a convex polyhedral cone in ${\Bbb R}^3$. The condition 
$g(e)=1$ defines a plane $H$ in ${\Bbb R}^3$ and the intersection 
$B=H \cap K$ is a {\it base} of $K$, that is, a polygon such that 
every $g \in K$ can be uniquely represented in the form 
$g=\lambda h$ for some $h \in B$. 

Our goal is to determine the structure of $K$.
 This is somewhat more complicated than in the 
2-dimensional situations of Sections 9-10. 
\proclaim{(11.2) Proposition} Let us define functions 
$$\split  &r_1={-np+n+p^2-2 \over n-2}, \\
          &r_2=1-2t, \\
          &r_3={2np-3p-2n-p^2-2t+6 \over n^2-5n+6}, \\
          &r_4={p+2t-2 \over n-2} \qquad \qquad \text{and} \\
          &r_{5o}={-2np+3p^2-3p+2tn+n-3 \over n^2-2n-3}.
\endsplit$$
Then 
\roster
\item If $e \in S_n$ is the identity, then
$$r_1(e)=r_2(e)=r_3(e)=r_4(e)=r_{5o}(e)=1;$$
\item If $n$ is even then $r_1, r_2, r_3$ and $r_4$ are the vertices 
(in consecutive order) of the planar quadrilateral
$B=\conv\bigl\{r_1, r_2, r_3, r_4 \bigr\}$ which is a base of the
central cone $K$;
\item If $n$ is odd then $r_1, r_2, r_3, r_4$ and $r_{5o}$ are 
the vertices (in consecutive order) of the planar pentagon
$B=\conv\bigl\{r_1, r_2, r_3, r_4, r_{5o} \bigr\}$ which is a base of 
the central cone $K$.
\endroster
\endproclaim
\demo{Proof} A function $g \in \spa\{g_1, g_2, g_3\}$ can be written as 
a linear combination $g=\alpha_1 g_1 + \alpha_2 g_2 + \alpha_3 g_3$.
Then $g(e)=\alpha_1(n-1) + \alpha_2 (n^2-2) -\alpha_3$ and the 
conditions $g(e) \geq g(\sigma)$ are written as 
$$\alpha_1(n-1) +\alpha_2(n^2-2) +\alpha_3 \geq 
\alpha_1\bigl(p(\sigma)-1\bigr) +\alpha_2 \bigl(p^2(\sigma)-2\bigr)+
\alpha_3 \bigl(1-2t(\sigma) \bigr),$$
which, for $\sigma \ne e$ are equivalent to
$$\alpha_1 +\alpha_2\bigl(n+p(\sigma)\bigr) +\alpha_3
{2t(\sigma) \over n-p(\sigma)} \geq 0.$$
Applying Lemma 9.2, we see that for even $n$, the system is equivalent to 
$$\split &\alpha_1 +n \alpha_2 \geq 0 \\
         &\alpha_1 + (2n-3) \alpha_2 \geq 0 \\
         &\alpha_1 + (2n-2) \alpha_2 + \alpha_3 \geq 0 \\ 
         &\alpha_1 + n \alpha_2 + \alpha_3 \geq 0 \endsplit \tag11.2.1$$
whereas for odd $n$, the system is equivalent to
$$\split &\alpha_1 +n \alpha_2 \geq 0 \\
         &\alpha_1 + (2n-3) \alpha_2 \geq 0 \\
         &\alpha_1 + (2n-2) \alpha_2 + \alpha_3 \geq 0 \\
         &\alpha_1 + (n+1) \alpha_2 + \alpha_3 \geq 0  \\
         &n\alpha_1 +n^2 \alpha_2 + (n-3) \alpha_3 \geq 0.
  \endsplit \tag11.2.2$$
The set of all feasible 3-tuples $(\alpha_1, \alpha_2, \alpha_3)$ is 
a polyhedral cone, which, for even $n$, has at most 4 extreme rays and 
for odd $n$ has at most 5 extreme rays. We call an inequality of 
(11.2.1)--(11.2.2) {\it active} on a particular tuple if it holds with 
equality.

It is readily verified that for even $n$ the following tuples span 
the extreme rays of the set of solutions to (11.2.1):
$$\split  \bigl(-n, \quad 1, \quad 0\bigr)  \qquad
            &\text{4th and 1st inequalities are active} \\
          \bigl(0, \quad 0, \quad 1 \bigr)  \qquad 
            &\text{1st and 2nd inequalities are active} \\
          \bigl(2n-3, \quad -1, \quad 1 \bigr) \qquad 
            &\text{2nd and 3d inequalities are active} \\ 
          \bigl(1, \quad 0, \quad -1 \bigr) \qquad 
            &\text{3d and 4th inequalities are active}
\endsplit$$
and that for odd $n$ the following tuples span the extreme 
rays of the set of solutions to (11.2.1):
$$\split \bigl(-n, \quad 1, \quad 0\bigr)  \qquad 
             &\text{5th and 1st inequalities are active} \\ 
         \bigl(0, \quad 0, \quad 1 \bigr)  \qquad 
            &\text{1st and 2nd inequalities are active} \\
          \bigl(2n-3, \quad -1, \quad 1 \bigr) \qquad 
            &\text{2nd and 3d inequalities are active} \\ 
           \bigl(1, \quad 0, \quad -1 \bigr) \qquad 
            &\text{3d and 4th inequalities are active} \\
            \bigl(-2n-3, \quad 3, \quad -n \bigr) \qquad
            &\text{4th and 5th inequalities are active}
\endsplit$$
We obtain $r_1,r_2,r_3,r_4$ and $r_{5o}$ by scaling the corresponding 
linear combinations $\alpha_1 g_1 +\alpha_2 g_2 +\alpha_3g_3$ so that 
the value at the identity is equal to 1 and hence $r_1,r_2, r_3, r_4$ and 
$r_{5o}$ lie on the same plane in $\spa\{g_1, g_2, g_3\}$.
{\hfill \hfill \hfill} \qed
\enddemo
\remark{(11.3) Remark}
One can observe that if $n$ is even then $r_{5o} \notin K$, for 
if $\sigma$ is a product of $n/2$ commuting transpositions,
so that $p(\sigma)=0$ and $t(\sigma)=n/2$, then 
$r_{5o}(\sigma)=(n^2+n-3)/(n^2-2n-3)>1=r_{5o}(e)$.
$$\epsffile{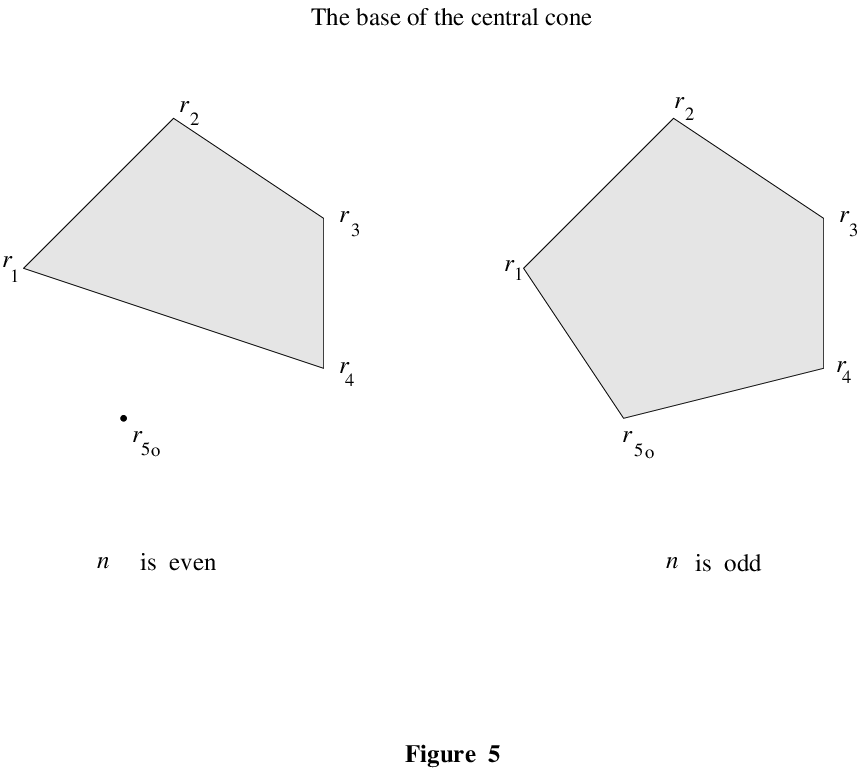}$$
We observe that function $r_3$ coincides with function $r_1$ of 
Lemma 10.2 (the symmetric QAP) and that function 
$r_2$ coincides with function $r_1$ of Lemma 9.3 
(the pure QAP). Function $r_4$ has 
a bullseye type distribution (see Section 2.2) whereas  
$r_1$ is a sharp spike (see Section 4.1). We have $r_1(\sigma)=1$ 
if and only if $\sigma=e$ or $\dist(e, \sigma)=n$ and $r_1(\sigma)<0$ 
for all other $\sigma$.
$$\epsffile{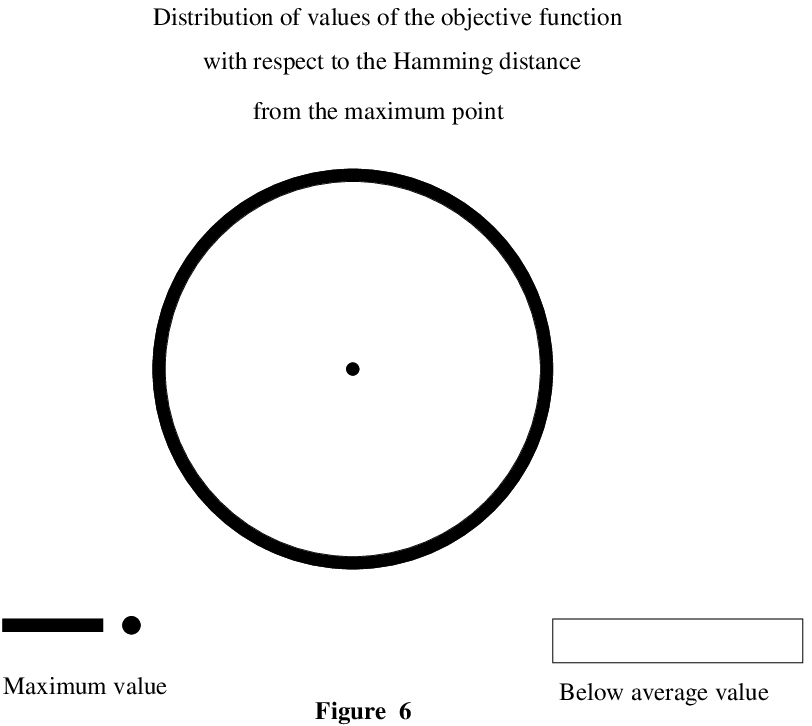}$$ 
Function $r_{5o}$ resembles a spike, but diluted.
\endremark    
Now we are getting ready to prove Theorem 5.1.
\proclaim{(11.4) Lemma} Let $g \in \spa\{g_1, g_2, g_3\}$
be a function such that $g(e)=1$.
For a $2 \leq k \leq n-2$, let $\sigma_k$ be a permutation such 
that $p(\sigma_k)=k$ and 
$t(\sigma_k)=0$, let $\eta$ be a permutation such that 
$p(\eta)=0$ and $t(\eta)=1$ and let $\theta$ be permutation such 
that $p(\theta)=t(\theta)=0$. Then 
$$\max\bigl\{g(\sigma_k), g(\eta), g(\theta) \bigr\} \geq 
{k-2 \over n^2-kn+k-2}.$$
\endproclaim
\demo{Proof} We can write 
$$g=\alpha_1 {p-1 \over n-1} + \alpha_2 {p^2-2 \over n^2-2} +\alpha_3(1-2t)$$
for some $\alpha_1, \alpha_2$ and $\alpha_3$ such that 
$\alpha_1 +\alpha_2 +\alpha_3=1$. 
Then 
$$\split &g(\sigma_k)=\alpha_1{k-1 \over n-1} +\alpha_2{k^2-2 \over n^2-2} +
\alpha_3 \\ 
&g(\eta)=-{\alpha_1 \over n-1} -\alpha_2 {2\over n^2-2} - \alpha_3 \\
&g(\theta)=-{\alpha_1 \over n-1} -\alpha_2{2 \over n^2-2} + \alpha_3. 
\endsplit$$
We observe that $g(\sigma_k)$, $g(\eta)$ and $g(\theta)$ are linear 
functions of $\alpha_1, \alpha_2$ and $\alpha_3$ and hence 
$$\ell(\alpha_1, \alpha_2, \alpha_3)=\max\bigl\{g(\sigma_k), g(\eta), 
g(\theta)\bigr\}$$ is a convex function on the 
plane $\alpha_1 +\alpha_2 +\alpha_3=1$. 

Moreover, for 
$$\alpha_1={k(1-n) \over n^2-nk+k-2}, \quad 
\alpha_2={n^2-2 \over n^2-nk+k-2} \quad \text{and} \quad 
\alpha_3=0 \tag11.4.1$$
we have 
$$g(\sigma_k)=g(\eta)=g(\theta)={k-2 \over n^2-nk+k-2}. \tag11.4.2$$
Let us prove that the minimum of $\ell(\alpha_1, \alpha_2, \alpha_3)$ on 
the plane $\alpha_1 +\alpha_2+\alpha_3=1$ is attained at (11.4.1).
Let 
$$\lambda_1={n^2-2n \over k^2-2k}, \quad \lambda_2={n^2-nk \over 2k-4} 
\quad \text{and} \quad \lambda_3={n^2-kn-2n+2k \over 2k}.$$
Then 
$$\lambda_1 g(\sigma_k) + \lambda_2 g(\eta) +\lambda_3 g(\theta)=
\alpha_1 +\alpha_2 +\alpha_3 =1 \quad \text{and} \quad 
\lambda_1,\lambda_2,\lambda_3>0.$$
Comparing this with (11.4.2), we conclude that there is no point 
$(\alpha_1, \alpha_2, \alpha_3)$ with 
$\alpha_1 +\alpha_2 +\alpha_3=1$ such that 
$$g(\sigma_k), g(\eta), g(\theta) < {k-2 \over n^2-nk+k-2}.$$ 
{\hfill \hfill \hfill} \qed
\enddemo
Now we are ready to prove Theorem 5.1.
\demo{Proof of Theorem 5.1} Without loss of generality, we may assume 
that the maximum value of $f_0$ is attained at the identity permutation
$e$. Excluding an obvious case of $f_0 \equiv 0$, by scaling $f$, if 
necessary, we may assume that $f_0(e)=1$. Let $g$ be the central projection 
of $f_0$. By Lemma 6.2, $g(e)=f_0(e)=1 \geq g(\sigma)$ for 
all $\sigma \in S_n$ and $\overline{g}=0$. By Proposition 7.5, $g$ 
must be a linear combination of the functions $\chi_n, \chi_{n-1,1},
\chi_{n-2,2}$ and $\chi_{n-2,1,1}$. Since $\overline{g}=0$, $g$ is 
a linear combination of non-trivial characters $\chi_{n-1,1}$,
$\chi_{n-2,2}$ and $\chi_{n-2,1,1}$ alone. Therefore, $g$ lies in the 
central cone: $g \in K$, see Definition 11.1.

Let
$X_k$ be the set of all permutations $\sigma$ such that $p(\sigma)=k$ and
$t(\sigma)=0$. As in the proof of Theorem 4.1, we conclude that 
$$|X_k| \geq {1 \over 5} {n! \over k!}.$$
Let $Y$ be the set of all permutations $\sigma$ such that 
$p(\sigma)=0$ and $t(\sigma)=1$. To choose a permutation $\sigma \in Y$,
one has to choose a transpositions in $\displaystyle 
{n \choose 2}$ ways and then an arbitrary permutation of 
the remaining $(n-2)$ symbols without fixed points and 2-cycles.
Using (6.4.2), we estimate
$$|Y| \geq {1 \over 5} {n! \over 2(n-2)!} (n-2)!={1\over 10} n!.$$ 
Let us choose a permutation $\sigma_k \in X_k$, a permutation 
$\eta \in Y$ and a permutation $\theta \in X_0$.
Let us choose $Z$ to be one of $X_k$, $X_0$ and $Y$, depending 
where the maximum value of $g(\sigma_k)$, $g(\eta)$ or $g(\theta)$ is 
attained. Hence
$$|Z| \geq {n! \over 5k!}.$$
The set $Z$ is a disjoint union of some conjugacy classes $X(\rho)$ and
for each $X(\rho)$ by (6.4.1) and Lemma 11.4, we have 
$$g(\rho)={1 \over |X(\rho)|} \sum_{\sigma \in X(\rho)} f_0(\sigma) 
\geq {k-2 \over n^2-kn+k-2}$$
and hence
$${1 \over |Z|} \sum_{\sigma \in Z} f_0(\sigma) \geq {k-2 \over n^2-kn+k-2}.$$
Applying Lemma 6.3 with $X=Z$ and $\beta=\beta(n,k)$, we 
conclude that 
$$\PP\Bigl\{\sigma \in S_n: f_0(\sigma) \geq \gamma\beta(n,k) \Bigr\}
\geq {(1-\gamma) \beta(n,k) \over 5k!}.$$
for all $n \geq 5$.
{\hfill \hfill \hfill} \qed
\enddemo

\head References \endhead

\item{1.} K. Anstreicher, N. Brixius, J.-P. Goux and J. Linderoth,
 Solving large quadratic assignment problems on computational grids,
{\it preprint}, 2000.

\item{2.} E. Arkin, R. Hassin and M. Sviridenko, 
Approximating the maximum quadratic assignment problem, 
{\it Inform. Process. Lett.}, {\bf 77} (2001), no. 1, 13--16.  

\item{3.} G. Ausiello, P. Crescenzi, G. Gambosi, V. Kann, 
A. Marchetti-Spaccamela, and M. Protasi, {\it Complexity and
Approximation. Combinatorial optimization problems and their
approximability properties},
Springer-Verlag, Berlin, 1999. 

\item{4.} 
A. Br{\" u}ngger, A. Marzetta, J. Clausen and M. Perregaard,
Solving large scale quadratic assignment problems in parallel with the
 search library ZRAM, {\it Journal of Parallel and Distributed Computing},
{\bf 50}, pp. 157-66, 1998.

\item{5.} 
R. Burkard, E. {\c C}ela, P. Pardalos and L. Pitsoulis,
The quadratic assignment problem, 
in: {\it Handbook of Combinatorial Optimization} 
(D.-Z. Du and 

P.M. Pardalos, eds.),
 Kluwer Academic Publishers, pp. 75-149, 1999.

\item{6.}
W. Fulton and J. Harris, 
{\it Representation Theory}, Springer-Verlag, New York, 1991.

\item{7.} I.P. Goulden and D.M. Jackson, Combinatorial
Enumeration, {\it Wiley-Interscience Series in Discrete Mathematics},
John Wiley $\&$ Sons, Inc., New York, 1983. 
\enddocument

\end